\documentclass[11pt,oneside,reqno]{amsart}
\usepackage{amssymb,epsfig,amsmath,amsthm,amsopn,amstext,amsfonts,amsbsy,graphicx}
\usepackage{color}
\usepackage[latin1]{inputenc}
\usepackage[T1]{fontenc}
\usepackage[english]{babel}
\input epsf

\def\Pt{\tilde p_t}
\def\Pi{\tilde p_\infty}

\newtheorem{theorem}{Theorem}
\newtheorem{proposition}{Proposition}
\newtheorem{lemma}{Lemma}

\newtheorem{corollary}{Corollary}

\title[Fluctuations of the front in a one dimensional model]{Fluctuations of the front in a one dimensional model
of $X+Y\to 2X$}

\author{Francis Comets$^{1,4}$, Jeremy Quastel$^{2}$ and Alejandro F. Ram\'\i rez$^{3,4}$}

\thanks{ AMS 2000 {\it subject classifications}. Primary  82C22, 82C41;
 secondary 82C24, 60K05, 60G50.}


\thanks{$^1$Partially supported by CNRS, UMR 7599.}

\thanks{$^2$Partially supported by NSERC, Canada}

\thanks{$^3$Partially supported by Fondo Nacional de Desarrollo Cient\'\i fico
y Tecnol\'ogico grant 1060738}

\thanks{$^4$Partially supported by ECOS-Conicyt grant CO5EO2}

\address[Francis Comets]{Laboratoire de Probabilit\'es et Mod\`eles
 Al\'eatoires\\
Universit\'e Paris 7- Denis Diderot\\
2, Place Jussieu\\
F-75 251 Paris Cedex 05, France}
\email{comets@math.jussieu.fr}

\address[Jeremy Quastel]{Departments of Mathematics and Statistics\\
University of Toronto\\
40 St. George Street\\
Toronto, Ontario M5S 1L2, Canada}
\email{quastel@math.toronto.edu}
\address[Alejandro F. Ram\'\i rez]{Facultad de Matem\'aticas\\
Pontificia Universidad Cat\'olica de Chile\\
Vicu\~na Mackenna 4860, Macul\\
Santiago, Chile}

\email{aramirez@mat.puc.cl}
\bigskip

\begin{document}

\begin{abstract}
We consider a model of the reaction $X+Y\to 2X$ on the
integer lattice in which $Y$ particles do not move while $X$ particles move as independent
continuous time, simple symmetric random walks.   $Y$ particles are 
 transformed instantaneously to $X$ particles upon contact.
We  start with a fixed number $a\ge 1$ of $Y$ particles at each site to the
 right of the origin, and define a class of  configurations of the $X$
  particles to the left of the origin having a finite $l^1$ norm
with a specified exponential weight. Starting from any configuration
of $X$ particles to the left of the origin within such a class, 
we prove a central limit
 theorem for the position of the rightmost visited site of the
 $X$ particles. 
 \end{abstract}

\maketitle

\section{Introduction}

We consider the following microscopic model of a combustive reaction
or epidemic
on the integer lattice $\mathbb Z$:
There are two types of particles; $X$ particles, which move as
independent,
continuous--time,  symmetric, nearest neighbor random walks of total jump rate
$2$; and
$Y$ particles which do not move.  Initially
the $Y$ particles occupy sites $1,2,\ldots$, with a fixed number $a\ge 1$
of $Y$ particles at each site.  Initially there is at least one $X$ particle
at $0$, and any  distribution of $X$ particles at sites $\ldots,-2,-1$
such that $\sum_{x\le 0}\eta (0,x)e^{\theta x}<\infty$, where
$\theta>0$ is a parameter that will be chosen small and
$\eta(0,x)$ is the number of $X$ particles at  $x\in{\mathbb Z}$
 at time $0$. 
   When an $X$ particle jumps to a site where there are
$Y$ particles, all $a$ of them immediately become $X$ particles and
start moving as rate $2$ continuous time symmetric random walks.
We are interested in the asymptotic behavior of the rightmost site
$r_t$ visited by the $X$ particles up to time $t$, which we call the {\it front}.  

Let $\eta(t,x)$ denote the number of $X$ particles at $x\in{\mathbb Z}$ at time $t\ge 0$.
Since there are always exactly $a$ of the $Y$ particles at each $x>r_t$
we do not have to keep track of them and we can just think of an $X$ particle as branching into
$a+1$ particles when it jumps to $r+1$, with the
result that there are $a+1$ particles
at the new rightmost visited site, $r+1$.
 A naive
 state space of our Markov process is 
\begin{equation}
\nonumber
 \mathbb S = \{ (r,\eta) ~:~r\in \mathbb Z, 
\eta\in {\mathbb N}^{\{\ldots, r-1,r\}}\},\end{equation}
with an infinitesimal generator acting over local functions  given by,
\begin{equation}
\nonumber
{\mathcal L}f(r,\eta) = \sum_{x,x+e \le r} \eta(x) (f(r,\eta-\delta_x + \delta_{x+e} ) -f(r,\eta))\end{equation}\begin{equation}
\nonumber + 
\eta(r) (f(r+1,\eta-\delta_r + (a+1)\delta_{r+1}) - f(r,\eta)).
\end{equation}
where $\delta_x$ denotes the configuration with one particle
at $x$.  Nevertheless, to avoid anomalies involving an explosion on the number
of particles per site, we will take as the state space
of our process,
\begin{equation}
\nonumber
\mathbb S'_\theta = \{ (r,\eta)\in \mathbb S ~:~\sum_{x\le r} e^{\theta(x-r)} \eta(x) <\infty\}.
\end{equation}
$\mathbb S'_\theta$ with, for example the metric 
$d((r,\eta),(r',\eta')) = |r-r'|+ \sum_{x\le 0} e^{\theta x} |\eta(x+r)-\eta'(x+r')|$,
is a Polish space.

We will show (see Section 2 and 6) that if initially
$(r,\eta)\in \mathbb S'_\theta$, with $r=0$ and $\eta(0,0)\ge 1$, then
 $(r_t,\eta(t))\in \mathbb S'_\theta$ and
furthermore  the process is Feller.  In \cite{rs} it is shown, 
for certain initial conditions, that
there exists  $v$, $0<v<\infty$,  such that a.s.,
\begin{equation}
\nonumber
\lim_{t\to\infty}{r_t}/{t}=v.
\end{equation}
We will give an
 alternate proof in dimension $d=1$ using the regeneration
time method (see Section 6) which works for arbitrary initial
data in $\mathbb S'_\theta$.  Note that this could also be proved using
the sub-additive ergodic theorem.

Our main results are:

\begin{theorem}
\label{theorem1}   (Central limit theorem) For $\theta>0$ small enough,
there exists $\sigma^2$ nonrandom,  $0<\sigma^2<\infty$, and independent of the
the initial conditions
$(0,\eta)\in{\mathbb S}'_\theta$, such that 
\begin{equation}
\label{process}
B_t^\epsilon:=\epsilon^{1/2}\left(r_{\epsilon^{-1}t}-\epsilon^{-1} vt\right),
\qquad\qquad t\ge 0,
\end{equation}
converges in law as $\epsilon\to 0$  to Brownian motion
with variance  $\sigma^2$.
\end{theorem}

\begin{theorem}
\label{theorem2}   (Ergodic theorem)  Consider the process as seen from the
front, $
\tau_{-r_t}\eta(t)
$.  For $\theta>0$ small enough,
 there exist exactly two invariant measures: One 
 supported on the configuration
with no particles, and another, $\mu_\infty$.  The domain of attraction of the first consists of exactly the configuration with no 
particles.  Any nontrivial configuration in $\mathbb S'_\theta$ is in 
the domain of the second; if we denote by $\mu_t$ the distribution 
of the process $
\tau_{-r_t}\eta(t)
$, then $\mu_t\to \mu_\infty$ in the sense of weak convergence of
probability measures.  
\end{theorem}

The model we are studying has been considered in the physics literature
(see \cite{panja} and references therein). Recently there has been a resurgence of interest
in such models because, especially in one and two dimensions, strong deviations
from mean field behavior were detected experimentally.

Mathematically much less is known.   \cite{betal} studies a model with at most one particle per site in which particles jump to neighboring sites at rate $\gamma/2$ and create particles at empty neighboring sites at rate $\frac12$.  Considering initial configurations with a rightmost particle, it is shown that 
viewed from the rightmost particle, the process has a unique invariant measure. Therefore the 
position of the rightmost particle grows linearly, in fact with a computable speed.

A discrete time version of our model is known in the probability
literature as the "frog model".   Shape theorems have been obtained
for the model on $\mathbb Z^d$ using methods based on the sub-additive
ergodic theorem (see \cite{amp2} and \cite{rs} for the continuous time
version and \cite{ampr} where the initial configuration of the $Y$
particles is random).  We prove the corresponding result for
arbitrary initial conditions in ${\mathbb S}'_\theta$ (see Section 6) which
could alternately be obtained with such methods.  However, using the
method of regeneration times we are able to obtain in addition the central
limit theorem for the position of the front and the ergodic
theorem for the law of the process as seen from the front.  The disadvantage of the
method is that it appears at the present time to be restricted to 
one dimensional systems.

In \cite{KS}, Kesten and Sidoravicius consider a model in which the $Y$ particles move as well. Let $D_X$ and $D_Y$ denote the
jump rates of the two types.  If $D_X=D_Y>0$ they prove a shape theorem in $\mathbb Z^d$.  
When $D_X\neq D_Y $ they can only obtain a linear upper bound. Note that \cite{mskb}  
observed experimentally that for one dimensional models of this type with exclusion, the speed does not depend
on $D_Y\ge 0$ but only on $D_X$ (as long $D_X>0$).

One of the aspects which makes these type of problems difficult is that the process as seen from the
front does not converge exponentially fast  to its equilibrium.  For example,
starting from one $X$  particle at the origin, the probability that the rightmost occupied site up to
time $t$ is still at the origin decays with $\mathcal O(t^{-1/2})$.
Hence, with such an initial condition, if $\mu_t$ is the law of the environment seen
from the the front at time $t$
and $\mu_\infty$ the (nontrivial) invariant measure of the process seen from the front,
 \begin{equation}
\nonumber
 || \mu_t - \mu_\infty ||_{TV} \ge \mathcal O(t^{-1/2}),\end{equation}
indicating that  we are in the gap-less case. In the physics literature such fronts are called {\it pulled fronts}
\cite{S}.

In \cite{cqr} we considered a preliminary model in which there was
a threshold:  Any particle which jumps to a site with $M$ particles
is immediately killed.  That model lacks the sub-additivity of the
present model.  On the other hand, it is considerably easier in
that case to define the renewal structure. The unboundedness in the particle configurations makes it particularly difficult to set up the renewal structure.  Essentially one has to show that at the regeneration time, one is
not in a bad situation in which there are an unusually large number of particles around.
Nevertheless, if uniform estimates on the
initial conditions are not obtained, then there is no finiteness of the first
and second moments of the corresponding regeneration times.
Therefore, to prove Theorem \ref{theorem1}, we have defined regeneration
times in terms of a modified renewal structure which provides a global
control on the number of particles per site far from the front.  
This difficulty in constructing regeneration times appears to be very 
common when dealing with dynamic environments (see, for example, \cite{bz}).

Regeneration time methods were already used by Kesten in \cite{kesten} to
study the invariant measure of an i.i.d. environment
 as seen from a one dimensional random walk  on that
 environment (RWRE).
Our approach to define the regeneration times in terms a
sequence of stopping times is inspired in the methods presented in
 \cite{sznitman-zerner} for multidimensional RWRE.
At a heuristic level,
regeneration occurs each time the front moves forward and the particles
behind it never catch it up later on. After such a time,
the behaviour of the front depends only on the $a$ newly created particles
sitting at the front at that time, but not on those behind the front at that time. The idea is to find an increasing sequence $\{\kappa_n:n\ge 1\}$
of regeneration times, having independent increments and 
such that the probability of the event $\{\kappa_n>t\}$ decreases
fast enough as $t\to\infty$ providing
 good enough integrability conditions. As in \cite{cqr},
in order to estimate the tails of the regeneration times, it is useful
to decouple particles initially on the front from those behind it.  Nevertheless, a crucial difficulty and difference in
the construction of the sequence of stopping times with respect to
\cite{cqr}, is that in this model the number of $X$ particles per site
is not bounded. This requires a control in terms of
some norm of the size of the cloud of particles behind the front.
To do so, we introduce at each time $t\ge 0$, an {\it exponential norm}
depending on the parameter $\theta$  and on an integer $z$,
which  is given by
$\sum_{x\le r_t}e^{\theta(x-r_t)}\eta_{z}(t,x)$. Here, $\eta_{z}(t,x)$
is the number of $X$ particles at site $x$ and at time $t$ which originated
from some branching (of an $X$ particle) at some site $y\le z$.
 This is a measure
of the magnitude of the density of particles from $r_t$ to $-\infty$, which
originated from some site $y\le z$.
We then  define a stopping time $S$ depending on an integer
length $L$, as the first hitting time
to a site of the form $r_0+jL$, $j\ge 1$, such that the exponential
norm of the particles originating to the left of $r_0+(j-1)L$
is small enough.  In \cite{cqr}, the corresponding stopping
time was defined simply as the first time the front advances $L$
steps to the right. One of the main difficulties of our
proof, is to show that the tails of the law of $S$
  provide good enough integrability conditions
for the corresponding regeneration times and the associated position of the front. We are able to do this
only for small  values of $\theta$ and large values of $L$:
we obtain polynomially decaying
tails of a  degree which increases linearly with $L$
 for the regeneration times $\{\kappa_n:n\ge 2\}$. It is conceivable that for a  fixed value
of $L$, the optimal bound for the corresponding regeneration times is
indeed of power law type (see \cite{sz} for a discussion of this problem 
within the context of transient multidimensional RWRE).

In the next section, we will define the notion of exponential norm, and
the labeled and auxiliary processes, which will be needed subsequently
to define the renewal structure. In Section 3, the renewal structure
is defined, following the algorithmic approach of \cite{sznitman-zerner}. Here it
is proved that the regeneration times, define sequences with independent
increments, and except for the first term, are identically distributed.
This is used in Section 4, to prove the law of large numbers, the
central limit theorem in Theorem \ref{theorem1}, and Theorem \ref{theorem2}.
In Section 5, the crucial estimates which ensure the finiteness of the
second moments of the i.i.d. sequences defined through the regeneration times
are derived. Of particular importance is Lemma \ref{fundamental}, which
shows that the tails of the stopping time $S$ are small enough. Finally, in Section 6, it is proved that the process
is Feller on ${\mathbb S}'_\theta$.  Note that in related models (see \cite{KS}) it is not known
whether the  Feller property holds. Throughout the paper a generic constant
will be denoted by $C$.

\smallskip
\section{Setup and preliminary definitions}


The process will be constructed out of a large collection of independent,
continuous time,  symmetric, simple random walks, each with jump
rate $2$.  For each site $x\le r$, we have a countable
collection of these: $\{Y_{x,1}, Y_{x,2},\ldots\}$. For each site $x>r$, we need only $a$ of them:
$\{Y_{x,1},\ldots,Y_{x,a}\}$. Assume that $Y_{x,i}(0)=x$.

First we construct the process for  finite initial conditions $(r,\eta)$
i.e., those in which
$\eta$ has only a finite number of particles.
 
 For each $x\le r$, and $i\le \eta(x)$,
 let $Z_{x,i} (t)=Y_{x,i}(t)$.
 Let $\tau_1$ be the first time that one of the random walks
$Z_{x,i}(t)$,   $x\le r$, hits  $r+1$. For
$0\le t<\tau_1$, let $r_t=r$ and $\eta(z,t)=\sum_{x\le r, ~i}  1 (Z_{x,i}(t) = z)$. 

At time $\tau_1$  we add
$a$ particles, $\{Z_{r+1,1},\ldots, Z_{r+1,a}\}$,
where $Z_{r+1,i}(t)=Y_{r+1,i}(t-\tau_1)$. 
Let $\tau_2$ be the first time  that one of the random walks
$Z_{x,i}(t)$,   $x\le r+1$, hits  $r+2$.
 For $\tau_1\le t<\tau_2$, 
let $r_t=r+1$ and $\eta(z,t)= \sum_{x\le r+1, ~i}  1( Z_{x,i}(t) = z)$.

 Continuing in this way, we define the process $\{(r_t,\eta(t):t\ge 0\}$ for finite initial
conditions and the sequence of stopping times $\{\tau_n:n\ge 1\}$.  In Section \ref{fellerp} we will show that the definition makes
sense.  In particular, one has to show that $\lim_{n\to\infty}
\tau_n=\infty$ with probability
one.

For general $(r,\eta)\in {\mathbb S}'_\theta$, with arbitrary $\theta$, we construct the process
by taking limits of approximations with finite initial conditions. For each $\ell=1,2,\ldots$,  let $\eta^\ell(x)=0$ if $x\le r-\ell$, and
$\eta^\ell(x)=\eta(x)$ if $r-\ell<x\le r$. 
 Consider the process
 $\{(r^\ell_t,\eta^\ell(t)):t\ge 0\}$ starting from this finite
 initial condition.  In Section \ref{fellerp} we will prove
 
 \begin{proposition}
\label{propfeller}
 For every $(r,\eta)\in{\mathbb S}'_\theta$ and
$t\ge 0$, $r_t=\lim_{\ell\to\infty} r^\ell_t$ 
 and $\eta(t,x) = \lim_{\ell\to\infty}\eta^\ell(t,x)$ exist, are finite
 a.s.  and $(r_t,\eta(t))\in {\mathbb S}'_\theta$. The limit is a Markov process with  Feller 
 semi-group $P_t f(r,\eta)= E_{r,\eta}[ f(r_t,\eta(t))]$ on $C(\mathbb S'_\theta)$, where $E_{r,\eta}$ is the expectation associated to the
joint law $P_{r,\eta}$ of $\{(r_t,\eta(t)):t\ge 0\}$.
 \end{proposition}
 
\subsection{Auxiliary process}\label{aux}
Let  
\begin{equation}
\label{mcond}
M= 4(a +5).
\end{equation}
Let now $r\in{\mathbb Z}$, define
$\nu_0:=0$ and $\nu_1$ as the first time one of the
random walks $\{Y_{r,i}:1\le i\le a\}$, hits the site $r+1$.
Next, define $\nu_2$ as the first time one of the random
walks $\{Y_{z,i}:r\le z\le r+1, 1\le i\le a\}$ hits the site $r+2$.
In general, for $k\ge 2$, we define $\nu_k$ as the first time
one of the  random walks
 $\{Y_{z,i}:r\lor (r+k-M)\le z\le r+k-1, 1\le i\le a\}$,
 hits the
site $r+k$. For $n\in {\mathbb N}$, let
\begin{equation}\nonumber
\tilde r_t^r:=r+n,\qquad{\rm if}\qquad 
\sum_{k=0}^n\nu_k\le t<\sum_{k=0}^{n+1}\nu_k.
\end{equation}
Now, observing that for each $1\le j\le M-1$, the random variables
$\{\nu_{Mk+j}:k\ge 1\}$ are independent and have finite moments since $M\ge 3$, we see that
a.s. (see also \cite{cqr}), 
\begin{equation}
\label{vel-alfa}
\lim_{t\to\infty} \tilde r^r_t/t=:\alpha>0.
\end{equation}

\subsection{Labeled process}
We enlarge the state space of the stochastic combustion
process so that particles carry labels indicating at which site
they originated. Each particle will have a starting
position $z\in{\mathbb Z}$ and label $(x,i), x\in{\mathbb Z},
i\in\{1,\ldots ,a\}$ describing its birthplace, allowing
the possibility that $z\ne x$. Throughout the sequel, we will
adopt the convention of calling $x$ the site where the
particle {\it originated}, whereas $z$ the site where
the particle was initially.

We fix at time $0$, an $r\in{\mathbb Z}$ representing the rightmost
visited site, and a subset ${\mathcal I}(0)$ of the labels
$(x,i)$ with $x\le r$, representing the set of labels of particles
at time $0$. To each one of these labels we assign a position
$z=Z_{x,i}(0)\le r$ which is the position at time $t=0$ of that
particle. The position at time $t$ is $Z_{x,i}(t)=Y_{x,i}(t)+z-x$.
Now, the first time a particle jumps to site $r+1$, the
labels $\{(r+1,1),\ldots ,(r+1,a)\}$ are
added to the set of labels of particles.
Let us call $\rho_1$ the time this happens. The trajectories $Z_{r+1,i}(t)$
 of these new particles are then equal to $Y_{r+1,i}(t-\rho_1)$
for $t\ge \rho_1$. Similarly, for $k\ge 2$, $\rho_1+\cdots+\rho_k$
will be the first time a particle jumps to $r+k$ adding at that
time the labels $\{(r+k,1),\ldots ,(r+k,a)\}$ to ${\mathcal I}$,
with trajectories $Z_{r+k,i}(t)=Y_{r+k,i}(t-\rho_1-\cdots\rho_k)$
for $t\ge\rho_1+\cdots \rho_k$.

Now denote by ${\mathcal I}(t)$ the set of labels of particles at time $t$
and by ${\mathcal Z}(t):=\{Z_{x,i}(t):(x,i)\in{\mathcal I}(t)\}$
their corresponding positions. We
assume that initially the set of labels of particles
includes at least one with $x=r$. Then, the rightmost visited
site is defined as $r_t=\sup\{x:(x,i)\in{\mathcal I}(t)\}$. 
Call ${\mathbb L}$ the triples $(r,{\mathcal I},{\mathcal Z})$
of integers $r$, labels ${\mathcal I}\subset
\{(x,i):x\le r, 1\le i\le a\}$ and position function
${\mathcal Z}:{\mathcal I}\to\{\ldots,r-2,r-1,r\}$. 
The unlabeled process defined in the previous section is
just the particle count 
\begin{equation}\nonumber
\eta(y,t)=\sum_{(x,i)\in{\mathcal I}(t)} 1(Z_{x,i}(t)=y).
\end{equation}
For $\theta>0$, let us now denote by ${\mathbb L}_\theta$ the set of triples
$(r,{\mathcal I},{\mathcal Z})\in{\mathbb L}$ such that
$(r,\eta)\in {\mathbb S}'_\theta$.
Then define

\begin{equation}\nonumber
\tilde{\mathbb S}_\theta:=\left\{(r,{\mathcal I},{\mathcal Z})
\subset{\mathbb L}_\theta:\max_{(x,i)\in{\mathcal I}}x=r\right\}.
\end{equation}

From Proposition \ref{propfeller}, note that if
$w_0=(r_0,{\mathcal I}(0),
{\mathcal Z}(0))\in \tilde{\mathbb S}_\theta$ then
$w_t=(r_t,{\mathcal I}(t),
{\mathcal Z}(t))\in \tilde{\mathbb S}_\theta$ for $t\ge 0$.
We now define the {\it labeled process} starting from $w_0$
as the triple
 $\{w_t:t\ge 0\}=\{(r_t,{\mathcal I}(t), {\mathcal Z}(t)):t\ge 0\}$, with a law
given by a probability measure ${\mathbb P}_w$ defined on the
Skorohod space $D([0,\infty);\tilde{\mathbb S}_\theta)$. 
Throughout this paper, we will occasionally use the
notation ${\mathbb P}_{r,eta}$ to denote any law ${\mathbb P}_w$
with an initial condition $w$ compatible with $r$ and the
particle count $\eta$.

Using sub-additivity we have the following result (see also Lemma 3 of \cite{cqr}),

\begin{lemma}
\label{subad} Suppose that
$(r,1),\ldots ,(r,a)\in{\mathcal I}(0)$, all initially at $r$.
Then $\rho_k\le \nu_k$.
\end{lemma}

Let us now define ${\mathcal R}(t)$ as the set of labels obtained
after removing from ${\mathcal I}(t)$ all labels $(x,i)$
with $x<r=\sup\{y:(y,i)\in{\mathcal I}(0)\}$. We define
for $y\le r_t$ the particle count
\begin{equation}\label{zetadef}
\zeta(t,y):=\sum_{(x,i)\in{\mathcal R}(t)}1(Z_{x,i}(t)=y).
\end{equation}

\subsection{Exponential density norm of particles}
Assume that
the initial condition of our process is $(r,\eta)$.
Let us also fix two integers $z_1,z_2$,  such that $z_1<z_2\le r-1$
 and follow the individual  particles which originated at  $z_1< y\le z_2$:
\begin{equation}\nonumber
\eta_{z_1,z_2}(t,y):= \sum_{(x,i):z_1<x\le z_2}1(Z_{x,i}(t)=y),
\end{equation}
 We will also write $\eta_z(t,y)$ for
$\eta_{-\infty,z}(t,y)$.
 We will use the notation,

\begin{equation}\nonumber
m_{z_1,z_2}(t):=\sum_{x=z_1+1}^{z_2}\eta_{z_1,z_2}(x,t),
\end{equation}
to denote the total number of such particles which are still in the same interval
at time $t$.

 For $\theta>0$ and $t\ge 0$
define,
\begin{equation}\nonumber
\phi_z(t,r,\eta):=\sum_{x\in{\mathbb Z}}e^{\theta (x-r_t)}\eta_z(t,x),
\end{equation}
which we will call the {\it exponential density norm} of particles.
Sometimes we will write $\phi_z(t)$ instead
of $\phi_z(t,r,\eta)$.

\section{The renewal structure }
Let us now define the renewal structure that will be used to
define the regeneration times. The exponential density norm of particles
will be an important ingredient and will enable us to control
the number of particles far from the front.
Let us now fix  some integer $L$ satisfying 
\begin{equation}
\label{lbound}
aL\ge M,
\end{equation}
and real numbers $\theta,\alpha_1$ and $\alpha_2$ satisfying
\begin{equation}
\label{teta}
0<2\sinh 2\theta < \alpha_1<\alpha_2 <\alpha=\lim_{t\to\infty} \tilde r^r_t/t .
\end{equation}
Let us now consider the labeled process $w_t$ with its natural
filtration ${\mathcal F}_t$ with an initial condition
$w_0\in\tilde{\mathbb S}_\theta$ having particles with labels
$(r,1),\ldots ,(r,a)$ at site $r$,   and
any allowable configuration of particles with labels to the left of $r$. 
Call $\eta(0)$ the initial particle count corresponding
to $w_0$. 

Define the stopping times,
\begin{equation}\nonumber
W:=\inf\{t\ge 0: \phi_{r-L} (t,r,\eta(0))\ge e^{\theta (\lfloor
\alpha_1 t\rfloor-(r_t-r))}\},
\end{equation}
and
\begin{equation}\nonumber
V:=\inf\{t\ge 0: \max_{r-L< z< r}\max_{1\le i\le a}
Z_{z,i}(t)>\lfloor\alpha_1 t\rfloor +r\}.
\end{equation}
When $W=\infty$,
none of the particles initially to the left of $r-L$ ever touches 
the line $\lfloor\alpha_1 t\rfloor+r$. Define,

\begin{equation}
\nonumber
U:=\inf \{t\ge  0: \tilde r^r_t -r< \lfloor \alpha_2 t\rfloor\}.
\end{equation}
We then let

\begin{equation}\nonumber
D:=\min\{U,V,W\}.
\end{equation}
We will  also need to define $U\circ\theta_s, V\circ\theta_s$ and $W\circ\theta_s$ as the first times $U,V$ or $W$ happen starting from the initial
condition $w_s$ for
$s\ge 0$, 
%
and $D\circ\theta_s:=\min\{U\circ\theta_s,V\circ\theta_s,W\circ\theta_s\}$.

For each $y\in{\mathbb Z}$, let

\begin{equation}\nonumber
T_y:=\inf\{t\ge 0:r_t\ge y\}.
\end{equation}
Fix $p$ such that 
\begin{equation}
\label{pbound}
0<pe^\theta<1.\end{equation}
We will furthermore impose the following additional condition
on $L$,
\begin{equation}
\label{lc}
(a-1)e^{-L\theta}<p
\end{equation}
 Now define for  $x\ge r$,
\begin{equation}\label{jay}
J_x:=\inf\{j\ge 1: \phi_{x+(j-1)L}(T_{x+jL})\le p
\ {\rm and}\ 
m_{x+(j-1)L,x+jL}(T_{x+jL})\ge aL/2\}.
\end{equation}
Define the sequence of ${\mathcal F}_t$-stopping times,
$\{S_k:k\ge 0\}$ and $\{D_k:k\ge 1\}$ as follows. Let $S_0:=0$ and $R_0=r$.
Then define

\begin{equation}\nonumber
S_1:=T_{R_0+J_{R_0}L} \qquad D_1:=D\circ\theta_{S_1}+S_1,
\qquad R_1:=r_{D_1}\end{equation}
and for $k\ge 1$,
\begin{equation} \nonumber S_{k+1}:=T_{R_k+J_{R_k}L}\qquad
 D_{k+1}:=D\circ\theta_{S_{k+1}}+S_{k+1},\qquad R_{k+1}=r_{D_{k+1}}
\end{equation}
Let
\begin{equation} \nonumber
K:=\inf\{k\ge 1:S_k<\infty ,D_k=\infty\},
\end{equation}
and define the {\it regeneration time} 
\begin{equation} 
\label{six}
\kappa:=S_{K},
\end{equation}
with the understanding that $\kappa=\infty$ on the event 
$\{k\ge 1:S_k<\infty ,D_k=\infty\}=\emptyset$.
Note that $\kappa$ is {\em not} a stopping time with
respect to ${\mathcal F}_t$. 

Define ${\mathcal G}$, the information up to time $\kappa$, defined
as the completion with respect to ${\mathbb P}_w$ of the
smallest $\sigma$-algebra containing all sets of the
form $\{\kappa\le t\}\cap A, A\in{\mathcal F}_t$. 
\begin{proposition}
\label{moment2}
For every initial condition $w\in\tilde{\mathbb S}_\theta$ with at least one
particle at the rightmost visited site, 

\begin{equation} 
\label{weight1}
 \kappa<\infty,\qquad \qquad {\mathbb P}_w-{\rm a.s.}
\end{equation}
Furthermore,
if $a\delta_0$ denotes a configuration with rightmost
visited site $0$ such that the number of particles
at $0$ is $a$ and the number of particles at each site $x<0$ is $0$,

\begin{equation}
\label{moments-delta}
{\mathbb E}_{a\delta_0}\left[ \kappa^2|U=\infty\right]<\infty\qquad
{ and} \qquad
{\mathbb E}_{a\delta_0}\left[ r^2_{\kappa}|U=\infty\right]<\infty.
\end{equation}

\end{proposition} 
\noindent Proposition \ref{moment2} will be proved in Subsection
 \ref{regen-tp}.
Recall the definition (\ref{zetadef}) of $\zeta$.
The key observation is 

\begin{proposition} \label{prop1} Let $A$ be a Borel subset of
$D([0,\infty);{\mathbb S}'_\theta)$ and $w\in\tilde{\mathbb S}_\theta$.
 Then,
\begin{equation} \nonumber
{\mathbb P}_w[\tau_{-r_\kappa}\zeta(\kappa+\cdot)\in A|
{\mathcal G}]={\mathbb P}_{a\delta_0}[\eta(\cdot )\in A|U=\infty].
\end{equation}
\end{proposition}

\begin{proof}[Proof] We have to show that for any
 $B\in{\mathcal G}$,
\begin{equation}\label{cc}
{\mathbb P}_{w}[B, \{\tau_{-r_{\kappa}} \zeta(\kappa+\cdot)\in A\} ]
={\mathbb P}_{w}[B]~{\mathbb P}_{a\delta_0}[ \eta(\cdot)\in A ~|~ U=\infty].
\end{equation}
Now, using (\ref{weight1}),
\begin{eqnarray}
\nonumber
&\ &{\mathbb P}_{w}[B, 
\{ \tau_{-r_{\kappa}} \zeta(\kappa+\cdot)\in A\} ]
={\mathbb P}_{w}[\{\kappa<\infty\}, B, \{ \tau_{-r_{\kappa}} 
\zeta(\kappa+\cdot)\in A\} ]
\\
\nonumber
&=&\sum_{k=1}^\infty {\mathbb P}_{w}\left[\{S_k<\infty, D_k=\infty\}, B,  
 \tau_{-r_{\kappa}}\zeta(\kappa+\cdot)\in A \right]
\\
\label{one}
&=&
\sum_{k=1}^\infty\sum_{x\in {\mathbb Z}}{\mathbb P}_{w}[r_{S_k}=x,
S_k<\infty, D_k=\infty , B,
  \tau_{-x} \zeta(S_k+\cdot)\in A ].
\end{eqnarray}
 From the definition of $\mathcal G$ there is an event $B_k\in {\mathcal F}_{S_k}$
such that
$B_k=B$ on $\kappa=S_k$.
Therefore, we can continue developing (\ref{one}) to obtain,
\begin{eqnarray}
\nonumber
&=&\sum_{k=1}^\infty\sum_{x\in {\mathbb Z}}
{\mathbb P}_{w}\left[r_{S_k}=x, S_k<\infty, D_k=\infty, B_k, \tau_{-x}
\zeta(S_k+\cdot)\in A \right]\\
&=&\sum_{k,x}
{\mathbb E}_{w}\left[ 1(r_{S_k}=x, S_k<\infty, B_k)
{\mathbb P}_{w} \left[ D_k=\infty, \tau_{-x}\zeta(S_k+\cdot)\in A
~|~{\mathcal F}_{S_k}\right] 
\right],
\end{eqnarray}
where ${\mathbb E}_w$ is the expectation defined by ${\mathbb P}_w$.
But on the events $S_k<\infty$ and $r_{S_k}=x$, we have that
\begin{equation}
\label{zeta-eta}
\zeta_{w}(S_k+\cdot )=\eta_{a\delta_{x}}(\cdot)
\end{equation}
when $U_k=V_k=W_k=\infty$, and that $\eta_{a\delta_{x}}(\cdot)$ is independent of the configuration
of particles initially to the left of $x$. Here, $a\delta_x$, is the configuration
  with rightmost visited site $x$ and with $a$ particles at site $x$ with labels $(x,1),\ldots,(x,a-1)$ and none elsewhere.
Indeed, on the event $V_k=W_k=\infty$,
 the particles with initial positions $z$ to the left of $x$,
 are never to the right of $\lfloor \alpha_1 t\rfloor +x$. 
And on the event $U_k=\infty$, 
the front $r_t$ is always to the right of $\lfloor\alpha_2t\rfloor+x$
and hence of $\lfloor\alpha_1t\rfloor+x$. Therefore, there is no
effect of the particles initially to the left of $x$ on the front $r_t$,
so that $\zeta_w(S_k+\cdot)=\eta_{a\delta_{x}}(\cdot)$.
Then, (\ref{zeta-eta}) combined with
 the independence of $U_k$ and $V_k\land W_k$ given ${\mathcal F}_{S_k}$,
the translation invariance,  and the strong Markov property imply that on the
events $S_k<\infty$ and $r_{S_k}=x$,
\begin{eqnarray}
\nonumber
 &&
{\mathbb P}_{w} \left[  U_k=\infty,V_k\land W_k=\infty,\tau_{-x}
\zeta(S_k+\cdot)\in A
~|~{\mathcal F}_{S_k}\right]\\
\nonumber
&=&{\mathbb P}_{w} \left[  U_k=\infty,\tau_{-x}
\eta_{a\delta_{x}}(\cdot)\in A ~|~ {\mathcal F}_{S_k}\right]
 {\mathbb P}_{w}\left[V_k\land W_k=\infty ~|~ {\mathcal F}_{S_k}\right]\\
& =&
{\mathbb P}_{a\delta_0}[ U=\infty,\eta(\cdot)\in A]{\mathbb P}_{w}[  V_k
\land W_k=\infty ~|~{\mathcal F}_{S_k}] . 
\end{eqnarray}
Summarizing, we have,
\begin{eqnarray}
\nonumber
&&
{\mathbb P}_{w}[B,\tau_{-r_{\kappa}} 
\zeta(\kappa+\cdot)\in A ]\\
\label{aa}
&&
=  {\mathbb P}_{a\delta_0}[ U=\infty,\eta(\cdot)\in A]
\sum_{k,x} 
{\mathbb P}_{w}[ V_k\land W_k=\infty,r_{S_k}=x,S_k<\infty, B_k].
\end{eqnarray}
Letting $A={\mathbb S}'_\theta$ gives 
\begin{equation}
\label{bb}
 {\mathbb P}_{w}[B]={\mathbb  P}_{a\delta_0}[U=\infty] \sum_{k,x} 
{\mathbb P}_{w}[ V_k\land W_k=\infty,r_{S_k}=x,S_k<\infty, B_k].
\end{equation}
(\ref{aa}) and (\ref{bb}) together imply (\ref{cc}).\end{proof}

  Now define  $\kappa_1\le\kappa_2\le\cdots $ by $\kappa_1:=\kappa$ and for $n\ge 1$
\begin{equation} \nonumber
\kappa_{n+1}:=\kappa_n + \kappa ( w_{\kappa_n+\cdot}).
\end{equation}
where $\kappa( w_{\kappa_n+\cdot})$ is the
regeneration time starting from $ w_{\kappa_n+\cdot}$ and
we set $\kappa_{n+1}=\infty$ on $\kappa_n=\infty$ for $n\ge 1$. We will call
$\kappa_1$ the {\it first regeneration time} and $\kappa_n$ the {\it $n$-th regeneration time}.

For each $n\ge 1$ we define the 
$\sigma$-algebra,
${\mathcal G}_n$,
as the completion with respect to ${\mathbb P}_{w}$ of the smallest $\sigma$-algebra containing all
 sets of the form $\{\kappa_1\le t_1\}\cap
\cdots\cap\{\kappa_n\le t_n\}\cap A$, $A\in{\mathcal F}_{t_n}$.
Now, noting that $\{\kappa_1=\infty\}$ is a null event for ${\mathbb P}_w$
one can see that $\{U<\infty\}\cap\{\kappa_1<\infty\}=\{\tilde r_U\le 
r_{\kappa_1}\}\cap\{\kappa_1<\infty\}\in{\mathcal G}_1$ (see
Lemma 5 of \cite{cqr}). Hence,
 $\{U=\infty\}\in{\mathcal G}_1$.
So we have the following general version of Proposition \ref{prop1},
\begin{proposition}
\label{cmp}
  Let $A$ be a Borel subset of 
$ D([0,\infty);{\mathbb S}'_\theta)$ and $w\in\tilde{\mathbb S}_\theta$. Then,
\begin{equation} \nonumber
{\mathbb P}_{w}[\tau_{-r_{\kappa_n}}\zeta(\kappa_n+\cdot)
\in A  ~|~{\mathcal G}_n]
={\mathbb P}_{a\delta_0}[  \eta(\cdot)\in A ~|~ U=\infty].
\end{equation}
\end{proposition}
We can now describe the full renewal structure.
\begin{corollary}
\label{iid}
 Let $w\in \tilde{\mathbb S}_\theta$.
(i)  Under ${\mathbb P}_w$, $\kappa_1,\kappa_2-\kappa_1, \kappa_3-\kappa_2, \ldots$ are independent, and $\kappa_2-\kappa_1, \kappa_3-\kappa_2, \ldots$ are identically distributed with law identical to that of $\kappa_1$ under
${\mathbb P}_{a\delta_0}[\cdot |U=\infty]$.
(ii) Under ${\mathbb P}_w$, $ r_{\cdot\land\kappa_1}, r_{(\kappa_1+\cdot)\land\kappa_2} - r_{\kappa_1}, r_{(\kappa_2+\cdot)\land\kappa_3} - r_{\kappa_2}, \ldots$ are independent, and $r_{(\kappa_1+\cdot)\land\kappa_2} - r_{\kappa_1}, r_{(\kappa_2+\cdot)\land\kappa_3} - r_{\kappa_2}, \ldots$ are identically distributed with law identical to that of $r_{\kappa_1}$ under
${\mathbb P}_{a\delta_0}[\cdot|U=\infty]$ .
\end{corollary}

\section{Limit theorems}
\label{limits}
We  now use the renewal structure to prove the law of large
numbers and the central limit theorem for $r_t$.
Throughout, we will consider an
 an initial condition $(0,\eta)\in {\mathbb S}'_\theta$
 such that $\eta(0,0)\ge 1$.


\subsection{Law of Large Numbers} 
We will prove that,

\begin{equation}
\label{as-lim}
\lim_{t\to\infty} \frac{r_t}{t}=v:=\frac{{\mathbb E}_{a\delta_0}[r_{\kappa_1}
|U=\infty]}{{\mathbb E}_{a\delta_0}[\kappa_1|U=\infty]}.
\end{equation}
Note that we have that $\kappa_1<\infty$,
 ${\mathbb P}_{0,\eta}$-a.s. Hence, by Corollary \ref{iid}
a.s.

\begin{equation}
\label{kappa}
\lim_{n\to\infty}\frac{\kappa_n}{n}={\mathbb E}_{a\delta_0}[\kappa_1
|U=\infty],\quad{\rm and}\quad
\lim_{n\to\infty}\frac{r_{\kappa_n}}{n}={\mathbb E}_{a\delta_0}[r_{\kappa_1}
|U=\infty].\end{equation}
Now, for $t\ge 0$, define $n_t:=\sup\{n\ge 0:\kappa_n\le t\}$, with the
convention $\kappa_0=0$. From (\ref{kappa}) we see that a.s. $n_t<\infty$.
Also, $\lim_{t\to\infty}r_{\kappa_{n_t}}/t=v$. The limit
(\ref{as-lim}) now follows from the observation,
\begin{equation} \nonumber
\lim_{t\to\infty} t^{-1}|r_t-r_{\kappa_{n_t}}|=0,
\end{equation}
which is a consequence of the inequality $|r_t-r_{\kappa_{n_t}}|
\le |r_{\kappa_{n_t+1}}-r_{\kappa_{n_t}}|$ and the fact that
$\lim_{t\to\infty}r_{\kappa_{n_t}}/t=v$ a.s.

\subsection{Central limit theorem}
 Consider the quantity
$B_t^\epsilon$ defined in (\ref{process}) 
and
\begin{eqnarray}
&\Sigma_m:=\sum_{j=1}^m R_j,
&\end{eqnarray}
where $R_j:=r_{\kappa_{j+1}}-r_{\kappa_j}-(\kappa_{j+1}-\kappa_j)v$.
Now, for
$0\le t\le T<\infty$,
\begin{eqnarray}
 \nonumber
&|B_t^\epsilon-\epsilon^{1/2}\Sigma_{n_{t/\epsilon}}|\\
&
\le 2\epsilon^{1/2}\sup_{0\le n\le n_{\lfloor \epsilon^{-1}T\rfloor}}(r_{\kappa_{n+1}}-r_{\kappa_n})
+
2v\epsilon^{1/2}\sup_{0\le n\le n_{\lfloor \epsilon^{-1}T\rfloor} }(\kappa_{n+1}-\kappa_n).
\end{eqnarray}
On the other hand, from Corollary \ref{moment2}, we can conclude that
for every $u>0$,
\begin{equation} \nonumber
\lim_{\epsilon\to 0}{\mathbb P}_{a\delta_0}[
\epsilon^{1/2}\sup_{0\le n\le n_{\lfloor \epsilon^{-1}T\rfloor} }(\kappa_{n+1}-\kappa_n)>u]=0.
\end{equation}
Hence, in probability
\begin{equation}
\label{error1}
\sup_{0\le n\le n_{\lfloor \epsilon^{-1}T\rfloor}}
\epsilon^{1/2}(\kappa_{n+1}-\kappa_n)\to 0.
\end{equation}
and 
\begin{equation} \nonumber
\sup_{0\le n\le n_{\lfloor \epsilon^{-1}T\rfloor}}
\epsilon^{1/2}(r_{\kappa_{n+1}}-r_{\kappa_n})\to 0.
\end{equation}
This proves that $B_t^\epsilon-\epsilon^{1/2}\Sigma_{n_{\epsilon^{-1}t}}$ converges to
$0$ in probability, uniformly on compact sets of $t$. From Donsker's invariance principle, we know that $\sqrt{\epsilon}\Sigma_{\cdot /\epsilon}$
converges in law to a Brownian motion with variance
 ${\mathbb E}_{a\delta_0}[(r_{\kappa_1}-\kappa_1v)^2|U=\infty]$,
where $\Sigma_s, s\ge 0,$ now stands for the linear interpolation of $\Sigma_m,m\ge 0$.  Using that
 $\lim_{t\to\infty} n_t/t =1/ {\mathbb E}_{a\delta_0}[\kappa_1|U=\infty]$
we can conclude that
as $\epsilon\to 0$, $B_t^\epsilon$ converges to a
Brownian motion with  variance,

\begin{equation}
\label{sigma}
\sigma^2:=
\frac{{\mathbb E}_{a\delta_0}[(r_{\kappa_1}-\kappa_1v)^2|U=\infty]}
{{\mathbb E}_{a\delta_0}[\kappa_1|U=\infty]}.
\end{equation}

\subsection{ Non-degeneracy of the variance} We will show that
$\sigma^2>0$. It is enough to
show that there exists some $\beta$, $0<\beta<v$ such that,

\begin{equation} \nonumber
{\mathbb P}_{a\delta_0}[r_{\kappa_1}=L,L\beta^{-1}\le\kappa_1|U=\infty]>0.
\end{equation}
Now,
\begin{equation} \nonumber 
{\mathbb P}_{a\delta_0}[r_{\kappa_1}=L,L\beta^{-1}\le \kappa_1, U=\infty]\ge
{\mathbb P}_{a\delta_0}[L\beta^{-1}<S_1< U, D\circ \theta_{S_1}=\infty].
\end{equation}
But the right hand side  can be written as
\begin{equation} \nonumber
{\mathbb E}_{a\delta_0}[1(L\beta^{-1}<S_1<U ){\mathbb E}_{a\delta_0}[  
1(\min\{V\circ\theta_{S_1},W\circ \theta_{S_1}\}=\infty)1( U\circ \theta_{S_1}=\infty) ~|~{\mathcal F}_{S_1}]].
\end{equation}
Now note that given ${\mathcal F}_{S_1}$, $U\circ\theta_{S_1}$,
 $V\circ\theta_{S_1}$ and $W\circ\theta_{S_1}$ are independent.
Hence,

\begin{eqnarray} &
\mathbb E_{a\delta_0}[  1(\min\{V\circ \theta_{S_1},W\circ \theta_{S_1}\}
=\infty)1( U\circ \theta_{S_1}=\infty) ~|~\mathcal F_{S_1}] \nonumber\\ &
\!\!\!\!=\mathbb P_{a\delta_0}[V\circ \theta_{S_1}
=\infty ~|~\mathcal F_{S_1}]
\mathbb P_{a\delta_0}[ W\circ \theta_{S_1}
=\infty ~|~\mathcal F_{S_1}]
\mathbb P_{a\delta_0}[   U\circ \theta_{S_1}=\infty ~|~\mathcal F_{S_1}].
\end{eqnarray}
This implies that,

\begin{equation} \nonumber
{\mathbb P}_{a\delta_0}[L\beta^{-1}<S_1< U, D\circ \theta_{S_1}=\infty]\\
\ge C{\mathbb P}_{a\delta_0}[L\beta^{-1}<S_1< U],
\end{equation}
for some constant $C>0$. Now, we have to show that
${\mathbb P}_{a\delta_0}[L\beta^{-1}<S_1< U]>0$.
Note that the event $\{L\beta^{-1}<S_1<U\}$ contains the
following event: one of the initial $a$ particles at $0$
jumps to site $1$ at some time $v_1$, such that
 $\beta^{-1}<v_1<2\beta^{-1}$; the other $a-1$ particles 
initially at $0$ stay at the same site during the time
interval $[0,2L\beta^{-1}]$; at time $v_1$, one of the
$a$ particles originating at site $1$ jumps to site $2$
at some time $v_2+v_1$ such that $\beta^{-1}<v_2<2\beta^{-1}$;
the other $a-1$ particles born at site $1$ stay at the
same site during the time interval $[0,2L\beta^{-1}]$;
in general, if $k$ is such that $3\le k\le L$, at
time
$v_k+v_{k-1}+\cdots +v_1$ one of the particles born
at site $k$ moves to site $k+1$, and $\beta^{-1}<v_k<2\beta^{-1}$;
all other $a-1$ particles born at site $k$ stay at the same site
during the time interval $[0,2L\beta^{-1}]$. Note that 
 $T_L=v_1+\cdots+v_L$ and at this time we have $\phi_0(T_L)
\le (a-1)e^{-L\theta}$. By (\ref{lc}) this quantity
is smaller than $p$. It is easy to see
that the above described event has positive probability.

\smallskip
\subsection{Ergodic theorem}

Let $\Pt$ be the law at time $t$ of the process as seen from the front
\begin{equation} \nonumber
\tau_{-r_t} \eta(t) \in  
\tilde\Omega: =  \{ 0,1,2,\ldots\}^{{\mathbb Z}_- }  
\end{equation}
under ${\mathbb P}_{0,\eta(0)}$. Note that $\tau_{-r_t} \eta(t)$ is itself
a Markov process with infinitesimal generator
\begin{equation} \nonumber
\tilde {\mathcal L} f(\eta)=
\eta(0) [f(\tau_{-1}(\eta - \delta_0)+a\delta_0)-f(\eta)]  +
\sum_{x , y \leq 0,\atop |x-y|=1}
\eta(x) [f(\eta-\delta_x+\delta_y)-f(\eta)] 
\end{equation}

Let $f$ be a  bounded continuous local function $f$ on $\tilde\Omega$
Denote by $\ell(f)$ the smallest integer $\ell$ such that 
$f(\eta)$ does not depend on $\eta(x), x < -\ell$. The formula
 \begin{equation} \label{def:Pi}
\int_{\Omega_0} f d \Pi =
\frac{ 
{\mathbb E}_{a\delta_0}[  \int_{\kappa_N}^{\kappa_{N+1}}f(\tau_{-r_s} \eta(s))ds
 ~|~ U=\infty]}
{{\mathbb E}_{a\delta_0}[ \kappa_1 ~|~ U=\infty]}\;,\quad
N(\alpha_2-\alpha_1)>\ell(f)
 \end{equation}
defines a probability measure $\Pi$ on $\tilde\Omega$. 
The righthand 
side of (\ref{def:Pi}) does not depend on $N$ provided that 
condition $N(\alpha_2-\alpha_1)>\ell(f)$ holds. This shows that the 
family of probability measures defined on finite cylinders by this formula 
is consistent. 

\begin{theorem}
\label{th:env}
We have $\Pt \to \Pi$ weakly as $t \to \infty$, and $\Pi$ is invariant
for $\tilde L$.
\end{theorem}

\begin{proof} Let $f$ be bounded and continuous on $\tilde \Omega$.
To prove convergence, first note
that the last term in the decomposition
\begin{equation} \nonumber
\int_{\Omega_0} f d \Pt= 
{\mathbb E}_{0,\eta(0)}[\kappa_{N+1} \leq t, f(\tau_{-r_t} \eta(t))]
+ 
{\mathbb E}_{0,\eta(0)}[\kappa_{N+1} > t, f(\tau_{-r_t} \eta(t))]
\end{equation}
vanishes as $t \to \infty$. Also,
\begin{eqnarray} 
\nonumber
&{\mathbb E}_{0,\eta(0)}[\kappa_{N+1} \leq t, f(\tau_{-r_t} \eta(t))]\\
\nonumber &=\sum_{k \geq 1, x \in {\mathbb Z}} 
{\mathbb E}_{0,\eta(0)}[\kappa_{N+k} \leq t< \kappa_{N+k+1}, r_{\kappa_k}=x,
f(\tau_{-r_t} \eta(t))]\\
\nonumber &=\sum_{k \geq 1, x \in {\mathbb Z}} 
{\mathbb E}_{0,\eta(0)}\Big[r_{\kappa_k}=x,\;
{\mathbb E}_{\eta(0)}[\kappa_{N+k} \leq t< \kappa_{N+k+1}, 
f(\tau_{-r_t} \eta(t)) | {\mathcal G}_k ]\Big]\\
 &\!\!\!\!\!\! =\! \sum_{k \geq 1, x \in {\mathbb Z}} 
{\mathbb E}_{0,\eta(0)}\Big[r_{\kappa_k}\!\! =\! x,\!\! \;
{\mathbb E}_{0,\eta(0)}\!\left[\kappa_{N+k} \leq t< \!\kappa_{N+k+1}, 
\!f\!\left(
\tau_{-r_t}
\zeta^{(k)}
(t-\kappa_k)\right) \!\! ~|~\!\! 
{\mathcal G}_k \right]\!\Big]
\end{eqnarray}
where $\zeta^{(k)}$ is a short notation for  $\zeta(\kappa_k+\cdot)$.
Note that we have used that $N(\alpha_2-\alpha_1)>\ell(f)$.
By Proposition   \ref{cmp}, this quantity is equal to
\begin{eqnarray}
\nonumber
\!\!\!\!\!\!\!\!\!\!\!\!\!\!\!\!\!\!\!\!\!\!\!\!\!\!\!\!\!\!\!\!\!\!\!\!
\sum_{k \geq 1, x \in {\mathbb Z}} \int_0^t 
{\mathbb P}_{0,\eta(0)}[r_{s}=x, \kappa_k\in ds] \qquad\qquad\qquad\qquad
\qquad\qquad\qquad\qquad\qquad\qquad\\
\nonumber
 \qquad\qquad\qquad\qquad \times
{\mathbb E}_{a\delta_0}\left[\kappa_{N} \leq t-s < \kappa_{N+1}, 
f\left(
\tau_{-r_{t-s}}
\eta
(t-s)\right) ~|~ 
U=\infty \right] 
\\
\nonumber
{\stackrel{u=t-s}{=}}
\sum_{k \geq 1, x \in {\mathbb Z}} \int_0^t 
{\mathbb P}_{0,\eta(0)}[r_{u}=x, t\!-\!\kappa_k\in du]
\qquad\qquad\qquad\qquad\qquad\qquad\\
 \nonumber
\qquad\qquad\qquad\qquad \times
{\mathbb E}_{a\delta_0}\left[\kappa_{N} \leq u < \kappa_{N+1}, 
f\left(
\tau_{-r_{u}}
\eta
(u)\right) ~|~ 
U=\infty \right]\\
\label{star}
= \int_0^t {\mathcal N}_t(du) F_f(u)
\end{eqnarray}
where 
\begin{equation} \nonumber
{\mathcal N}_t([0,u])= 
\sum_{k \geq 1}  
{\mathbb P}_{0,\eta(0)}[ \kappa_k\in [t-u,t]]
\end{equation}
and 
\begin{equation} \nonumber
 F_f(u)={\mathbb E}_{a\delta_0}\left[\kappa_{N} \leq u < \kappa_{N+1}, 
f\left(
\tau_{-r_{u}}
\eta
(u)\right) ~|~ 
U=\infty \right].
\end{equation}
We will use the following  renewal theorem (Theorem 6.2 in \cite{thor}).  To state the
theorem we say that a random walk 
\begin{equation} \nonumber
S_n= S_0+X_1+\cdots+X_n, \qquad n=0,1,2,\ldots
\end{equation}
i.e. $X_1,X_2,\ldots$ are i.i.d. and independent of $S_0$, is a renewal process if $S_0$ is nonnegative and $X_k$ are strictly positive.
We say it 
has {\em spread out} step-lengths if there exists an $r\ge 1$ and
a nonnegative measurable function $m$ such that $\int_{\mathbb R} m(x)dx>0$ and
\begin{equation} \nonumber
P(X_1+\cdots+X_r\in A)\ge \int_A m(x) dx,
\end{equation}
for all Borel sets $B$.
\begin{theorem} {\it ( Renewal theorem)}. \label{thor}
Let $S$ be a renewal process with spread out step lengths and  $E[X_1]<\infty$.
For Borel sets $B$, let
\begin{equation} \nonumber
N(B)=\sum_{k=0}^\infty \mathbf 1_{\{S_k\in B\}}.
\end{equation}
Then for each $h\in [0,\infty)$,
\begin{equation} \nonumber
E[ N(t+B)] \to |B| / E[ X_1]
\end{equation}
uniformly over Borel sets $B\subset [0,h]$.  Here $|B|$ is the Lebesgue measure
of $B$.
\end{theorem} One can check the spread-out assumption in
 Theorem \ref{thor} as follows: With $T_L$ the time of $L$-th jump
for  the particle
with label $(0,1)$
first jumps, $A$ the event that all these $L$ jumps are to the right, $B$
the event that no other particle moves between times $0$ and $1$, we have  
for $0 <s<t<1$,
\begin{eqnarray}
\nonumber
{\mathbb P}_{0,\eta(0)}[\kappa_2-\kappa_1 \in (s,t]]&=&
{\mathbb P}_{a\delta_0}[ \kappa_1  \in (s,t] ~|~U=\infty]\\
\nonumber
&\geq & \frac
{{\mathbb P}_{a\delta_0} [T_L  \in (s,t], A, B, U \circ \theta_{1}=\infty,
V \circ \theta_{1}=\infty
]}{
{\mathbb P}_{a\delta_0}[U=\infty]} \\
&=&
C \int_s^t f_L(u) du,
\end{eqnarray}
with $f_L$ the $L$-fold convolution of the exponential density  with rate
$2$ and $C$ is a constant that we can check using independence satisfies
$C>0$. This shows that $\kappa_2-\kappa_1$ is
 spread-out.

Hence from the renewal theorem,
\begin{equation} \label{eq:renouv}
{\mathcal N}_t(B) \to {|B|}/
{{\mathbb E}_{a\delta_0}[ \kappa_1 ~|~ U=\infty]} \quad {\rm as} \; t \to 
\infty
\end{equation}
uniformly over Borel sets $B$ in any finite interval. 

Since $F_f(u)$ is bounded and measurable, we have from (\ref{star})
\begin{equation} \nonumber
\int_{\tilde\Omega} f d \Pt \to \int_{\Omega_0} f d \Pi
\end{equation}
Because the process is  Feller (Proposition \ref{propfeller}),  any limit measure is
invariant.
\end{proof}

\section{Expectations and variances of the regeneration times}
\label{expect}
\subsection{Bounds on $W$}

\begin{lemma}
\label{max}
 Let $\{X_t:t\ge 0 \}$ be a simple symmetric continuous time
rate $2$ random walk on ${\mathbb Z}$, such that $X_0=x$. Let
$M_t:=x+\sup_{0\le s\le t}\left|X_s-x\right|$. Then, for $t\ge 0$,

\begin{equation} \nonumber
E\left[ e^{\theta M_t}\right]\le 3e^{\theta x+2(\cosh\theta -1)t},
\end{equation}
where $E$ is the expectation defined by the law of the random walk.
\end{lemma}

\begin{proof} The reflection principle tells us that for every
integer $n\ge 0$,
$P[M_t\ge n]=2P[X_t>n ]+P[X_t=n]$. Hence
we have, $P[M_t=n]\le P[X_t=n]+2P[X_t=n+1]$. Therefore,
$E[e^{\theta M_t}]\le E[e^{\theta X_t}]+2e^{-\theta}E[e^{\theta X_t}]
\le 3E[e^{\theta X_t}]$. Finally remark that
$E[e^{\theta X_t}]=e^{\theta x+2(\cosh\theta -1)t}$.

\end{proof}

We will in several occasions consider the
random process,

\begin{equation} \nonumber
M_{x,i}(t):=Z_{x,i}(0)+\sup_{0\le s\le t}|Z_{x,i}(s)-Z_{x,i}(0)|,
\end{equation}
defined for each $(x,i)\in{\mathcal I}$.
Furthermore, we will need to define for each initial
condition $(r_0,{\mathcal I}(0),{\mathcal Z}(0))$ compatible
with a particle count $\eta(0)$ and each $z\le r_0-1$ the quantity,

\begin{equation} \nonumber
\psi_z(t,r_0,\eta(0)):=\sum_{(x,i)\in {\mathcal I}(0),x\le z}
e^{\theta (M_{x,i}(t)-r_t)}.
\end{equation}
Usually we will drop the argument, writing
$\psi(t)$ instead of $\psi(t,r_0,\eta(0))$. Let us
also note that since,

\begin{equation} \nonumber
\phi_z(t,r_0,\eta(0)):=\sum_{(x,i)\in {\mathcal I}(0),x\le z}
e^{\theta (Z_{x,i}(t)-r_t)}.
\end{equation}
it is true that,

\begin{equation}
\label{phi-psi}
\phi_z(t)\le \psi_z(t),
\end{equation}
for every $t\ge 0$ and $z\le r_0-1$.
Due to condition (\ref{teta}), and
the intermediate value theorem it is true that,

\begin{equation}
\label{intermediate}
\mu:=\theta\alpha_1-2(\cosh\theta -1)>0.
\end{equation}
This enables us to obtain the following exponential bound.

\begin{lemma}
\label{f}
 For all initial conditions $(r,\eta)$ such that
$\phi_{r-L}(0,r,\eta)<\infty$ and $t\ge 0$ we have that,

\begin{equation} \nonumber
{\mathbb P}_{r,\eta}\left[t< W<\infty\right]\le 
C\phi_{r-L}(0,r,\eta)
\exp\left\{-\mu t\right\},
\end{equation}
where $C=3e^\theta e^{2(\cosh\theta -1)}
 \frac{e^{\mu}}{1-e^{-\mu}}$.
\end{lemma}

\begin{proof} Without loss of generality we assume $r=0$.
Let us first note that,

\begin{equation} \nonumber
{\mathbb P}_w\left[t< W<\infty\right]
\le {\mathbb P}_w\left[\cup_{s\ge t}\left\{\phi_{-L}(s)\ge
e^{\theta(\lfloor\alpha_1 s\rfloor-r_s)}\right\}\right].
\end{equation}
From inequality (\ref{phi-psi}) and
the fact that $M_{x,i}(t)$ is nondecreasing in $t$, it follows using
Lemma \ref{max} that,

\begin{eqnarray}
\nonumber
{\mathbb P}_w\left[t< W<\infty\right]
&\le&\sum_{n=[t]}^{\infty} {\mathbb P}_w\left[
\sum_{(x,i)\in{\mathcal I}(0), x\le -L}e^{\theta M_{x,i}(n+1)}\ge
 e^{\theta\lfloor\alpha_1 n\rfloor}\right]\\
\nonumber
&\le &3\sum_{n=[t]}^{\infty}
e^{2(\cosh \theta -1)(n+1)-\theta\lfloor\alpha_1 n\rfloor}
\sum_{(x,i)\in{\mathcal I}(0),x\le -L} e^{\theta  Z_{x,i}(0)}\\
&\le &3\phi_{-L}(0)\sum_{n=[t]}^{\infty}
e^{2(\cosh \theta -1)(n+1)-\theta\lfloor\alpha_1 n\rfloor},
\end{eqnarray}
 Summing up
the last expression over $n$ we finish the
proof of the Lemma.
\end{proof}

 Define for $t\ge 0$,  and $z\le r$,

\begin{equation} \nonumber
N_z(t):=e^{\theta r_t-2(\cosh \theta -1)t}\phi_{z}(t).
\end{equation}

\begin{lemma}
\label{lemma4}  Consider an initial condition $(0,\eta)$
and an integer $z$ such that
$z\le 0$ and  $\phi_z(0)<\infty$. Then,
  $\{N_z(t):t\ge 0\}$ is an ${\mathcal F}_t$-martingale.
\end{lemma}

\begin{proof} Let us remark that,

\begin{equation} \nonumber
N_z(t)=\sum_{(x,i)\in{\mathcal I}(0), x\le z} e^{\theta 
Z_{x,i}(t)-2(\cosh\theta -1)t}.
\end{equation}
Now, each one of the terms in the above sum is an ${\mathcal F}_t$-martingale.
Furthermore, since $\phi_z(0)<\infty$, the martingales
 $\sum_{(x,i)\in{\mathcal I}(0), -n\le x\le z} e^{\theta 
Z_{x,i}(t)-2(\cosh\theta -1)t}$, converge in $L^1({\mathbb P}_w)$ norm
to $N_z(t)$ as $n\to\infty$. Thus, $\{N_z(t):t\ge 0\}$ is an
${\mathcal F}_t$-martingale.
\end{proof}
\medskip

\begin{lemma}
\label{uniform-f} 
 There is a $\delta>0$ such that for all initial
conditions $w$ with particle count $\eta$, initial position
of the front $r=0$, such that $\phi_{-L}(0,0,\eta)\le p$,

\begin{equation} \nonumber
 {\mathbb P}_{w}\left[W<\infty\right]< 1-\delta.
\end{equation}
\end{lemma}
\begin{proof} By inequality (\ref{intermediate}), note that,

\begin{equation}
\label{unif-1}
{\mathbb P}_{w}\left[W<\infty\right]
\le
{\mathbb E}_{w}\left[e^{\left(\theta\alpha_1-2(\cosh\theta -1)
\right)W}1(W<\infty)\right].
\end{equation}
Now,  from the definition of the exponential
density norm and of the stopping time $W$, the a.s. right-continuity of the
trajectories of the random walks, and Fatou's Lemma, it follows that
$e^{\theta(\lfloor\alpha_1 W\rfloor-r_{W})}\le
\phi_{-L}(W)$. Hence, from inequality (\ref{unif-1}) we
conclude that ${\mathbb P}_{w}\left[W<\infty\right]$
is bounded by,

\begin{equation}
e^\theta {\mathbb E}_{w}\left[e^{\theta r_W
-2(\cosh\theta -1)W}\phi_{-L}(W)1(W<\infty)\right]
=e^\theta {\mathbb E}_{w}[N_{-L}(W)1(W<\infty)].
\end{equation}
Now, note that
$ {\mathbb E}[N_{-L}(W)1(W<n)]\le {\mathbb E}_w[N_{-L}(n\land W)]$.
 Thus, by
the optional stopping  theorem and Fatou's Lemma,

\begin{equation} \nonumber
{\mathbb E}_{w}[N_{-L}(W)1(W<\infty)]\le \lim_{n\to\infty}{\mathbb E}_w
[N_{-L}(n\land W)]=N_{-L}(0)\le p.
\end{equation}
This and the condition $pe^\theta<1$, shows that ${\mathbb P}_w[W<\infty]<1$.
\end{proof}

\subsection{Bounds on $V$}

\begin{lemma}
\label{v}
 There is a $C$, $0<C<\infty$, such that for all initial conditions $w$
and all $t\ge 0$

\begin{equation} \nonumber {\mathbb P}_w\left[ t<V<\infty\right]\le
 C\exp\left\{-tC\right\}.
\end{equation}
\end{lemma}

\begin{proof} Without loss of generality
we assume that initially  $r=0$.
Note that the probability $P_w\left[t<V<\infty\right]$ is bounded
by the probability that one of the random walks born at a site
between $-L$ and $-1$ is at the right of
$\lfloor\alpha_1 s\rfloor$ at some time $s\ge t$. 
Now this probability is bounded by the worst case in which initially all
these random walks,
$aL$, are at site $0$. But this has probability,

\begin{equation} \nonumber
aL P[t<\tau<\infty],
\end{equation}
where $\tau:=\inf\{t\ge 0: X_t> \lfloor \alpha_1 t\rfloor\}$,
$\{X_t:t\ge 0\}$ is a continuous time simple symmetric
random walk on ${\mathbb Z}$, of total jump rate $2$, starting from $0$,
and $P$ is its law. It is easy to prove that
this probability is bounded by $C\exp\left\{-Ct\right\}$ for
some constant $C<\infty$ (for example, see Lemma 8 of \cite{cqr}).
\end{proof}

\begin{lemma}
\label{uniform-v} There is a $\delta>0$ such that
for all initial conditions $w$,

\begin{equation} \nonumber
{\mathbb P}_{w}\left[V<\infty\right]< 1-\delta.
\end{equation}
\end{lemma}
\begin{proof} Without loss of generality
we can assume that $r=0$. Note that the probability ${\mathbb P}_w[V<\infty]$
is upper bounded by the probability that a random walk within a group
of $aL$ independent ones all initially at site $x=0$, at some
time $t\ge 0$ is at the right of $\lfloor\alpha_1 t\rfloor$.
 But this probability is $1-\gamma^{aL}$,
where $\gamma$ is the probability that a single random walk
starting form $x=0$ never is at the right of the curve 
$\{\lfloor\alpha_1 t\rfloor:t\ge 0\}$.
By Lemma 8 of \cite{cqr} we know that $\gamma<1$.
\end{proof}

\subsection{Bounds on $U$} The following two lemmas can be proved
observing that at each instant of time $t\ge \nu_j$, with $j\ge M+1$,
 the auxiliary process
has at least $M\ge 20$ particles behind the front (see also
\cite{cqr}).

\begin{lemma}
\label{u} There is a constant $C$,  $0<C<\infty$, such that for all initial
conditions $w$ with particles $(r,1),\ldots , (r,a)$ at 
 the rightmost site $r$,
  and all
$t>0$

\begin{equation} \nonumber
 {\mathbb P}_{w}\left[ t<U<\infty\right]\le
Ct^{-M/2}.
\end{equation}
\end{lemma}

\begin{lemma}
\label{uniform-u} 
 There is a $\delta>0$ such that for all initial conditions
$w$ with particles $(r,1),\ldots , (r,a)$ at 
 the rightmost site $r$,

\begin{equation} \nonumber
 {\mathbb P}_{w}\left[U<\infty\right]< 1-\delta.
\end{equation}
\end{lemma}

\subsection{Bounds on $D$}
\label{boundsd}

The following lemma is elementary.

\begin{lemma}
\label{first} There is a constant $C$, $0<C<\infty$, such that 
for every $t>0$

\begin{equation} \nonumber
{\mathbb P}_{a\delta_0}[ \nu_1>t]\le Ct^{-a/2},
\end{equation}
while for every $j\ge M+1$ and $t>0$

\begin{equation} \nonumber
{\mathbb P}_{a\delta_0}[ \nu_j>t]\le Ct^{-M/2},
\end{equation}
so that ${\mathbb E}_{a\delta_0}[\nu_j^{M/2}]<\infty$.
\end{lemma} From here we obtain the following estimate.

\begin{lemma}
\label{slowdown}
 Let  $\beta$  be such that
$0<\beta<\alpha$.
Then there is a constant $C$,  $0<C<\infty$, such that
the following statements are true.

\begin{itemize}

\item[a)] Assume that $\eta$ has at least $a$ particles at $0$. Then,

\begin{equation} \nonumber
{\mathbb P}_{0,\eta}
\left[ T_n>  n/\beta  \right]\le Cn^{-a/2}.
\end{equation}

\item[b)] Assume that $\eta$ is such that $m_{-L,0}(0)\ge aL/2$. Then,

\begin{equation} \nonumber
{\mathbb P}_{0,\eta}
\left[ T_n>  n/\beta  \right]\le Cn^{-M/4}.
\end{equation}

\item[c)] Assume that $\eta$ is a configuration with at
least one particle. Then, for all $k\ge M$ we have,

\begin{equation} \nonumber
{\mathbb P}_{0,\eta}
\left[ T_{n+k}-T_k>  n/\beta  \right]\le Cn^{-M/4}.
\end{equation}

\end{itemize}
\end{lemma}

\begin{proof} 
Let us prove part $(a)$. First
remark that ${\mathbb P}_{0,\eta} [T_n>n/\beta]\le {\mathbb P}_{a\delta_0}[T_n>n/\beta]$.
Now, $T_n=\sum_{i=1}^n\rho_i$. 
Hence, by Lemma \ref{subad} we have $T_{n}\le
\sum_{j=1}^n\nu_j$. Therefore,
${\mathbb P}_{a\delta_0} \left[ T_{n}>  n/\beta
\right]\le 
{\mathbb P}_{a\delta_0}\left[\sum_{i=1}^{n} \nu_{i}>  n/\beta\right]$.

Choose now $\beta'$ such that $\beta<\beta'<\alpha$. Then since
$1/\beta=(1/\beta-1/\beta')+1/\beta'$ and $\nu_1$ is stochastically
larger than $\nu_j$ for $j\ge 2$, we have for $n\ge M+1$,

\begin{equation}
\label{ll1}
{\mathbb P}_{a\delta_0} \left[ T_{n}>  n/\beta\right]
\le
M{\mathbb P}_{a\delta_0}\left[ \nu_1> \frac{n}{M}\left(
 \frac{1}{\beta}-\frac{1}{\beta'}\right)\right]
+
{\mathbb P}_{a\delta_0}\left[\frac{1}{n}\sum_{i=M+1}^n \nu_i>
\frac{1}{\beta'}\right].
\end{equation}
But, ${\mathbb P}_{a\delta_0}\left[\frac{1}{n}\sum_{i=M+1}^n \nu_i>
\frac{1}{\beta'}\right]\le {\mathbb P}_{a\delta_0}\left[
\frac{1}{n}\sum_{i=M+1}^n \gamma_i>c
\right]$, where $\gamma_j:=\nu_j-1/\alpha$
and $c:=\frac{1}{\beta'}-\frac{1}{\alpha}>0$.
On the other hand, for each $0\le i<l$, $l=\lfloor (M+1)/(a+1)\rfloor$,
the random variables $\{\gamma_{kl+i}:k\ge 1\}$ are independent.
Thus,

\begin{equation} \nonumber
{\mathbb P}_{a\delta_0}\left[\frac{1}{n}\sum_{i=M+1}^n \nu_i>
\frac{1}{\beta'}\right]\le 
\sum_{i=0}^{l-1}{\mathbb P}_{a\delta_0}\left[
\frac{1}{n}\sum_{k:(M+1-i)/l\le k\le n} \gamma_{kl+i}>(a+1)c/M\right].
\end{equation}
Now,  for $q\ge 2$, if $X_1,X_2,\ldots$ are independent and
identically distributed random variables with mean zero, and
if $E[|X_i|^q]<\infty$, then $E[|\sum_{i=1}^n X_i|^q]
\le C n^{q/2}$ for some $C<\infty$ (see item 16, page 60 of \cite{petrov}).
Hence, since by Lemma \ref{first}
we have ${\mathbb E}_{a\delta_0}[\gamma_j^{M/2}]<\infty$,
it follows that the last expression of the above display is
bounded by, $Cn^{-M/4}$, for some other constant $C<\infty$. Finally
observe that $M\ge 2(a+1)$, and use again Lemma \ref{first}
to bound the first term of inequality (\ref{ll1})  to finish the proof.
The proofs of parts $(b)$ and $(c)$ are
similar using the inequality (\ref{lbound}) satisfied by the
parameters $M$ and $L$.
\end{proof}

Let us now obtain the estimates for the stopping time $D$.  From lemmas  \ref{f}, \ref{v} and \ref{u} we obtain
\begin{corollary}
\label{est-D-gen} There is a constant $C=C(p)$, $0<C<\infty$, such
that for all initial conditions $w$ with
$\phi_{-L}(0,w)\le p$, and with particles
$(r,1)\ldots ,(r,a)$ at the rightmost visited site $r$,  and for all $t> 0$,

\begin{equation} \nonumber
{\mathbb P}_{w}\left[ t<D<\infty
\right]\le Ct^{-M/2}.
\end{equation}
\end{corollary}
\medskip
We also have the following lemma.

\begin{lemma}
\label{uniform}
There is a $\delta>0$ such that,
for all initial
conditions $w$ with particle count $\eta$ and initial position
of the front $r=0$ such that $\phi_y(0,0,\eta)\le p$ and
with particles $(0,1),\ldots ,(0,a)$ at $0$,

\begin{equation} \nonumber
{\mathbb P}_{w}\left[D<\infty\right]< 1-\delta.
\end{equation}
\end{lemma}

\begin{proof} Since $W,V$ and $U$ are independent, 

\begin{equation}
{\mathbb P}_{w}\left[D<\infty\right]=1-
{\mathbb P}_{w}\left[W =\infty\right]
{\mathbb P}_{w}\left[V=\infty\right]
{\mathbb P}_{w}\left[U=\infty\right].
\end{equation}
Applying lemmas \ref{uniform-f}, \ref{uniform-v} and \ref{uniform-u},
we end up the proof.
\end{proof}

We finish this subsection with three lemmas and
a corollary which will be subsequently
used to obtain estimates for the stopping time $S$. The following
lemma will be proved in Section 6. 

\begin{lemma} 
\label{front1}
There are  constants $C$ and $\gamma_0$, $0<C<\infty$
 and  $\gamma_0>0$, such that for 
 all $w\in\tilde{\mathbb S}_\theta$, $\gamma\ge \gamma_0$ and $t\ge 0$,
\begin{equation}
\label{estimater}
{\mathbb P}_w[r_t\ge\gamma t]\le \phi_0(0,w) e^{-Ct}.
\end{equation}
\end{lemma}
\smallskip

\begin{lemma}
\label{rd} There is a constant $C=C(p)$, $0<C<\infty$,
such that for all initial conditions $w$ such that
$\phi_{-L}(0,w)\le p$,  and $t> 0$,

\begin{equation} \nonumber
{\mathbb P}_w\left[ r_D>  t,
D<\infty\right]\le Ct^{-M/2}.
\end{equation}
\end{lemma}
\begin{proof} Note that,

\begin{eqnarray}
\nonumber
 &{\mathbb P}_w\left[ r_{D}>\gamma t, D<\infty\right]
\le {\mathbb P}_w\left[ r_{D}>\gamma t,   D\le t \right]
+{\mathbb P}_w\left[ t<D<\infty\right]\\
&\le  
 {\mathbb P}_w\left[ r_t>\gamma t\right]
+{\mathbb P}_w\left[  t<D<\infty\right].
\end{eqnarray}
The statement now follows from 
(\ref{estimater}) of Lemma \ref{front1}, the
fact that $\phi_0(0,w)\le \phi_{-L}(0,w)+aL$ and Corollary \ref{est-D-gen}.

\end{proof}

\begin{lemma} 
\label{impbound}
Consider an initial condition $w$ with
rightmost visited site $r=0$,  at least $a$ particles at $0$
and such that $\phi_{-L}(0,w)\le p$. Then, ${\mathbb P}_w$-a.s. on the event $\{D<\infty\}$ we have,

\begin{equation} \nonumber
\phi_{-L}(D)\le e^{\theta}.
\end{equation}

\end{lemma}
\begin{proof} First note that by the assumption
$\phi_{-L}(0,w)\le p<1$, necessarily we have $D>0$. Now, by definition of $U$, note that whenever $t\le U<\infty$,
we have $\tilde r_t\ge \lfloor\alpha_2 t\rfloor$.
By Lemma \ref{subad} we have $r_t\ge \tilde r_t$. It follows that
$ r_t\ge  \lfloor\alpha_2 t\rfloor$. Therefore, if $t\le U<\infty$, we have

\begin{equation}
\label{last3}
\lfloor\alpha_1 t\rfloor-r_t\le -(\lfloor\alpha_2 t\rfloor-\lfloor \alpha_1 t\rfloor)\le 0.
\end{equation}
Therefore, if $D=U$, inequality (\ref{last3}) shows that
$\lfloor\alpha_1 D\rfloor-r_D\le 0$. Hence, since in this case with
probability one $D<W$, it follows that,

\begin{equation}
\label{last}
\phi_{-L}(D)\le e^{\theta(\lfloor\alpha_1D\rfloor -r_D)}\le  1.
\end{equation}
Similarly, if $D=V$, since $D<U$ and $D<W$ happen with
probability one, inequality (\ref{last})
still holds a.s. On the other hand, if $W<\infty$ we have,

\begin{equation}
\label{last2}
\phi_{-L}(W)\le e^\theta e^{\theta(\alpha_1 W-r_W)},
\end{equation}
since in the worst case scenario, at time $W$ all particles
jump one step to the right. Hence, if $D=W$, since with
probability one we have $D<U$, by inequality (\ref{last3}), the
exponent in the right hand side of  (\ref{last2}) is non-positive,
so that $\phi_{-L}(W)\le e^\theta$.
\end{proof}

\begin{corollary}
\label{psid} There is a constant $C$, $0<C<\infty$, such that
for all initial condition $w$ with rightmost visited site
$r=0$, such that $\phi_{-L}(0,w)\le p$ and at least the particles
with labels $(0,1),\ldots , (0,a)$ at $0$,

\begin{equation} \nonumber
{\mathbb E}_w[\phi_{r_D}(D), D<\infty]<C.
\end{equation}

\end{corollary}
\begin{proof} Placing ourselves in the worst case scenario were
all particles born between sites $-L$ and $r_D$ are at site $r_D$ at time $D$,
we see  that $\phi_{r_D}(D)\le \phi_{-L}(D)+a(L+r_D)$. Hence, by
Lemma \ref{impbound}, $\phi_{r_D}(D)\le e^{\theta}+a(L+r_D)$. Lemma
\ref{rd} together with the fact that $M\ge 3$ finishes the proof.

\end{proof}
\medskip

\subsection{Bounds on $S$}
We will now perform some key estimates which will let us 
obtain fast enough decay estimates for the tail probabilities
of $J$ in Lemma \ref{fundamental}.

\begin{lemma} 
\label{nol} There exists 
a constant $C$, $0<C<\infty$,
such that the following statements are satisfied.

\begin{itemize}

\item[a)] For all
 initial conditions $w\in\tilde{\mathbb S}_\theta$
 with at least $a$ 
particles at  the rightmost visited site, and all $n\ge 1$,

\begin{equation} \nonumber
{\mathbb P}_w[m_{0,n}(T_n)< an/2]\le C\frac{1}{n^{a/2}}. 
\end{equation}

\item[b)] For all initial conditions $w\in\tilde{\mathbb S}_\theta$
with at least $a$ particles at the rightmost visited site, and
all $n\ge 1$,

\begin{equation}
\label{extra2}
{\mathbb P}_w[ m_{r_D,r_D+n}(T_{r_D+n})<an/2]\le C\frac{1}{n^{aM/(4+2M)}}.
\end{equation}

\item[c)] For all nontrivial initial conditions $w\in
\tilde{\mathbb S}_\theta$, 

$$
\lim_{n\to\infty}{\mathbb P}_w[m_{-n,L}(T_L)<aL/2]=0.
$$
\end{itemize}
\end{lemma}
\begin{proof}[Proof of part $(a)$]
 Choose $0<\beta<\alpha$. Then,

\begin{equation}
\label{nn2}
{\mathbb P}_w\left[m_{0,n}(T_n)< \frac{an}{2}\right]\le 
{\mathbb P}_w
\left[m_{0,n}(T_n)<\frac{a n}{2}, T_n\le\frac{1}{\beta}n\right]
+{\mathbb P}_w\left[T_n>\frac{1}{\beta}n\right].
\end{equation}
Note that the event $\{m_{0,n}(T_n)<an/2, T_n\le n/\beta\}$ is
contained in the event that at least one particle
born at any of the sites $\lfloor n/2\rfloor,\lfloor n/2\rfloor+1,\ldots,n$ 
hits some site $x\le 0$ in
a time shorter than or equal to $n/\beta$. Hence, we can conclude that,

\begin{equation}
\label{nn}
{\mathbb P}_w
\left[m_{0,n}(T_n)< \frac{an}{2}, T_{0,n}\le\frac{1}{\beta}n\right]
\le a(n-\lfloor n/2\rfloor) P[M'_{n/\beta}\ge n/2],
\end{equation}
where $P$ is the law of a simple symmetric rate $2$
 random walk $\{X_t:t\ge 0\}$
on ${\mathbb Z}$ starting from $0$ and $M'_t:=sup_{0\le s\le t}X_s$.
Now, by the reflection principle, $P[M'_t\ge x]\le 2P[X_t\ge x]$.
Hence, from inequality (\ref{nn}), we see that
${\mathbb P}_w\left[m_{0,n}(T_n)<a n/2, T_{0,n}\le\frac{1}{\beta}n\right]$
is bounded by $a(n+1)P[X_{n/\beta}\ge n/2]$. But, for every $t\ge 0$ and
positive integer $x$, $P[X_t\ge x]\le e^{-2t I(x/(2t))}$, where
$I(u)=u\sinh^{-1}u-\sqrt{1+u^2}+1$. Hence,
$a(n+1)P[X_{n/\beta}\ge n/2]\le (a+1)(n+1)\exp \left\{-\frac{2n}{\beta} I(\beta/4)\right\}$.
Finally, using the inequality ${\mathbb P}_w[T_n>n/\beta]
\le{\mathbb P}_{a\delta_0}
[T_n>n/\beta]$, part $(a)$ of Lemma \ref{slowdown} to bound the second term of inequality
(\ref{nn2}) and using the fact that $(a+1)(n+1)\exp \left\{
-\frac{2n}{\beta} I(\beta/4)\right\}\le C/n^{a/2}$ for $n$ large enough,
we conclude the proof.

\smallskip

\noindent {\it Proof of part (b)}. By part (a) and Lemma \ref{rd},
${\mathbb P}_w[m_{r_D,r_D+n}(T_{r_D+n})<an/2]$ is upper bounded by,

\begin{eqnarray}
\nonumber
&\sum_{k:1\le k\le n}{\mathbb P}_w[m_{k,k+n}(T_{k+n})<an/2]
+{\mathbb P}_w[r_D>m,D<\infty]\\
\nonumber
&\le Cm\frac{1}{n^{a/2}}+C\frac{1}{m^{M/2}},
\end{eqnarray}
for some constant $C>0$ and for every $m\ge 1$. Choosing
$m=n^{\frac{a}{2+a}}$ we obtain (\ref{extra2}).

\smallskip

\noindent {\it Proof of part $(c)$}. Note that,

\begin{equation}
\label{extra5}
{\mathbb P}_w\left[m_{-n,L}(T_L)< \frac{aL}{2}\right]\le 
{\mathbb P}_w
\left[m_{-n,L}(T_L)<\frac{a L}{2}, T_L\le n\right]
+{\mathbb P}_w\left[T_L>n\right].
\end{equation}
Clearly
 $\lim_{n\to\infty}
{\mathbb P}_w\left[T_L>n\right]=0$. On the other hand, an
argument similar to the one used to derive (\ref{nn}),
shows that the first term of the righthand side of
(\ref{extra5}) is bounded by $aL P[M'_{n}\ge n]$, which
tends to $0$ as $n$ tends to $\infty$.
\end{proof}

Throughout the sequel, to simplify notation, we will
define on the event $\{D<\infty\}$ for each $n\ge 1$,

\begin{equation} \nonumber
F_n:=T_{r_D+Ln}-D.
\end{equation}

\begin{lemma}
\label{est-f}
For every $0<\beta<\alpha$, there exists a constant $C<\infty$
 depending only on $\beta$, such that for
all initial conditions $w$ with rightmost visited site $r$,
with at least $aL/2$ particles at a distance strictly smaller than $L$ to $r$,
 and such that $\phi_{r-L}(0,w)\le p$, and
for all natural $n \ge 1$,

\begin{equation} \nonumber
{\mathbb P}_{w}\left[ F_n>  
\frac{1}{\beta}L 
n, D<\infty \right]
\le C\frac{1}{(nL)^{M/4-1}}.
\end{equation}

\end{lemma}

\begin{proof} Without loss of generality
we can assume that initially $r=0$. Note that
${\mathbb P}_{0,\eta}
\left[ F_n> \frac{1}{\beta} Ln, D<\infty\right]$
 is upper-bounded by,

\begin{equation}
\label{f1}
\sum_{k:1\le k\le  Ln}
 {\mathbb P}_{w}\!\!\left[ F_n> \frac{1}{\beta}
Ln,r_D=k, D<\infty \right]
+{\mathbb P}_{w}\left[ r_D> Ln, D<\infty \right].
\end{equation}
Now, on the event $\{D<\infty\}$ we have that
$T_{r_D}\le D$ so that $F_n\le T_{r_D+Ln}-T_{r_D}$.
Hence,

\begin{equation} \nonumber
{\mathbb  P}_{w}\left[ F_n> 
\frac{1}{\beta} Ln,r_D=k, D<\infty \right]
\le 
 {\mathbb P}_{w}\left[ T_{k+Ln}-T_{k}>\frac{1}{\beta} Ln\right].
\end{equation}
Now, by part $(c)$ Lemma \ref{slowdown}, for all $k> M$ we have
 $ {\mathbb P}_{w}\left[ T_{k+Ln}-T_k>\frac{1}{\beta} Ln\right]\le\frac{C}{(nL)^{M/4}}$, for
some constant $C<\infty$. On the other
hand for $1\le k\le M$, 
${\mathbb P}_{w}\left[ T_{k+Ln}-T_k>\frac{1}{\beta} Ln\right]\le
{\mathbb P}_{w}\left[T_{M+Ln}>\frac{1}{\beta} Ln\right]$.
Thus,  by part $(b)$ of Lemma \ref{slowdown}, since the initial condition
$w$ has at least $aL/2$ particles to the right of $r=0$ at a distance strictly 
smaller
than $L$ to the origin, we know that
${\mathbb P}_{w}\left[T_{M+Ln}>\frac{1}{\beta} Ln\right]\le \frac{C}{(nL)^{M/4}}$,
for some other constant $C<\infty$.
 We therefore conclude that,

\begin{equation}
\label{f2}
\sum_{k:1\le k\le \gamma Ln}
 {\mathbb P}_{w}\!\!\left[ F_n> \frac{1}{\beta}
Ln,r_D=k, D<\infty \right]\le  C\frac{1}{(nL)^{M/4-1}}.
\end{equation}
Using Lemma \ref{rd} to estimate the second term of display (\ref{f1}) and
combining this with inequality (\ref{f2}) we finish the proof.
\end{proof}

Now we will be concerned with proving that given 
$D_{k-1}<\infty$, the
stopping time $S_k$ happens almost surely 
and has tails that decay fast enough.

\begin{lemma}
\label{series} Let $q\ge 1$ be an integer.
Consider a sequence $\{a_k:k\ge 1\}$ of non-negative
real numbers such that $\sum_{k=1}^\infty a_k<1$ and
$\sum_{k=1}^\infty k^qa_k<\infty$. Assume that $\{c_m:m\ge 1\}$
is a sequence such that,

\begin{equation}
\label{rec-1}
c_1\le a_1,
\end{equation}
and  for every $m\ge 2$ we have that,

\begin{equation}
\label{recurrence}
c_m\le a_{m}+\sum_{k=1}^{m-1}a_{m-k}c_k.
\end{equation}
Then,

\begin{equation} \nonumber
\sum_{k=1}^\infty k^q c_k<\infty.
\end{equation}
\end{lemma}

\begin{proof}  We will use induction on $0\le q'\le q$ to prove the lemma.
We introduce the notation $A_{q'}:=\sum_{k=1}^\infty k^{q'}a_k$
and $C_{q'}:=\sum_{k=1}^\infty k^{q'}c_k$.
Let us first show that if $A_0<\infty$ then $C_0<\infty$.
Let $n\ge 2$
 be a fixed natural. Summing up inequality (\ref{rec-1}) with
inequalities
 (\ref{recurrence}) from $m=2$ to $m=n$ we see that,

\begin{equation} \nonumber
\sum_{k=1}^{n-1} c_k (1-\sum_{j=1}^{n-k} a_j)+c_n\le A_0,
\end{equation}
Taking the limit when $n\to \infty$ above and  using 
 Fatou's Lemma we conclude that,

\begin{equation} \nonumber
\sum_{k=1}^\infty c_k\le\frac{A_0}{1-A_0}<\infty.
\end{equation}
Now assume that $C_{q'-1}<\infty$ for some $1\le q'\le q$. We will show
that then $C_{q'}<\infty$. 
Summing up inequality (\ref{rec-1}) with
inequalities
 (\ref{recurrence}), multiplied by $m^{q'}$, from $m=2$ to $m=n$ we see that,

\begin{equation}
\label{kkk}
\sum_{m=1}^n m^{q'} c_m \le A_{q'}+\sum_{m=2}^n\sum_{k=1}^{m-1} m^{q'} a_{m-k} c_k.
\end{equation}
Substituting the binomial expansion $m^{q'}=\sum_{i=0}^{q'}\binom{q'}{i} (m-k)^i
k^{q'-i}$ on (\ref{kkk}) and interchanging the order of the summations
on $m$ and on $k$, we conclude that,

\begin{equation} \nonumber
\sum_{m=2}^{n-1} m^{q'}c_m \left(1-\sum_{j=1}^{n-m}a_j\right)+n^q c_n
\le A_{q'}+\sum_{i=1}^{q'}\binom{q'}{i}\sum_{k=1}^{n-1} k^{q-1}c_k
\sum_{m=1}^{n-k}m^ia_m.
\end{equation}
Taking the limit when $n\to\infty$ and using Fatou's Lemma, we get

\begin{equation} \nonumber
C_{q'}\le\frac{A_{q'}+\sum_{i=1}^{q'}\binom{q'}{i} C_{q'-i} A_i}{1-A_0}<\infty.
\end{equation}
\end{proof}
\medskip

\begin{lemma} 
\label{psi-estimates}
Consider an initial condition  $w\in\tilde{\mathbb S}_\theta$
such that the rightmost visited site associated to $w$ is
$r=0$ and at least one particle.

\begin{itemize}

\item[a)] For every $h>0, s>0$ and $n\ge 1$
we have

 \begin{equation}
\label{psi-1}
{\mathbb P}_w\left[\psi_0(T_n)>
h, T_n< s\right]
\le 3\frac{\psi_0(0,w)}{ h} e^{2(\cosh\theta-1)s-\theta n}.
\end{equation}

\item[b)]  For every $h>0, s>0$, $k\ge 1$ and $n\ge k$
we have
\begin{equation}
\label{psi-2}
{\mathbb P}_w\left[\left.\psi_{k}(T_{n})-\psi_{k-L}(T_{n})>
h, T_{n}-T_{k}< s\right|{\mathcal F}_{T_{k}}\right]
\le 3\frac{aL}{h}
e^{2(\cosh\theta-1)s-\theta (n-k)}.
\end{equation}
\end{itemize}
\begin{proof} Note that of the event $T_n<s$, it is true
that $e^{\theta M_{x,i}(T_n)}\le e^{\theta M_{x,i}(s)}$. Therefore,
since 
$\psi_0(T_n)=e^{-\theta n }\sum_{(x,i),x\le 0}e^{\theta M_{x,i}(T_n)}$,
we have

\begin{equation} \nonumber
\psi_0(T_n)\le e^{-\theta n }\sum_{(x,i),x\le 0}e^{\theta M_{x,i}(s)}.
\end{equation}
Now, by Lemma \ref{max}  we have that
${\mathbb E}_w\left[\sum_{(x,i),x\le 0}e^{\theta M_{x,i}(s)}\right]
\le 3\psi_0(0,w)e^{2(\cosh\theta -1)s}$. Hence,

\begin{equation} \nonumber
{\mathbb E}_w\left[\psi_0(T_n)\right]\le 3\psi_0(0,w)
e^{2(\cosh\theta -1)s-\theta n }.
\end{equation}
Using Tchebyshev's inequality we obtain (\ref{psi-1}). A similar
argument, using the fact that $\psi_{k}(T_k)-\psi_{k-L}(T_k)\le 
aL$ proves
(\ref{psi-2}).
\end{proof}

\end{lemma}

\smallskip

We end up this section with the following result providing
a tail estimate for the law of $J_{r_D}$ (with $J_x$ for
$x$ integer, defined in (\ref{jay})). An important idea
in the proof is that essentially, the
event that the exponential norm $\phi_{nL}$ is larger
 than $p$, is contained on the event that some of 
 the exponential norms $\phi_{jL}$ are
 larger than $p/2^j$ for some  $0\le j\le n-1$.
\begin{lemma}
\label{fundamental}Assume that   $M$ and $L$ satisfy (\ref{mcond}) and (\ref{lbound}),
$\theta, \alpha_1,\alpha_2$ (\ref{teta}) 
and $p$ satisfies (\ref{pbound}).
Then,  there is a constant $C$, $0<C<\infty$, and an integer $L_0$
such that if $L\ge L_0$, the following statements are satisfied.

\begin{itemize}

\item[a)]
  Consider an initial condition $w\in\tilde{\mathbb S}_\theta$
with rightmost visited site $r=0$, such that the number
of live particles
at $0$ is $a$, and   such that $m_{-L,0}(0)\ge aL/2$.  Then,
for every $t> 0$,

\begin{equation} \nonumber
{\mathbb P}_w[J_{r_D}>t,D<\infty]<Ct^{3-M/4} .
\end{equation}
\item[b)]  Consider an initial condition $w\in\tilde{\mathbb S}_\theta$
with rightmost visited site $r=0$. Then, for every $t> 0$,

\begin{equation} \nonumber
{\mathbb P}_w[J_0 >t,U=\infty]<Ct^{3-M/4} .
\end{equation}

\end{itemize}
Furthermore, for every nontrivial initial condition
$w\in\tilde{\mathbb S}_\theta$,

\begin{equation}
\label{partc}
J_0<\infty,\qquad\qquad {\mathbb P}_w-{\rm a.s.}
\end{equation}

\end{lemma}
\begin{proof}[Proof of part $(a)$] Call  $F'_i:=T_{r_D+iL}=F_i+D$
for $i\ge 1$.
 For  $n= 1,2,\ldots$,
\begin{equation}
\label{tail1}
{\mathbb P}[J>n, D<\infty]
\le {\mathbb P}\left[ B_n, D<\infty\right],
\end{equation}
where we have dropped the subscripts on $\mathbb P_w$ and $J_{r_D}$ and defined
\begin{eqnarray}&
\nonumber
B_n:=\cap_{i=1}^n \left\{\psi_{r_D+(i-1)L}(F_i',w_D)>p\right\}\cup B'_i,&\\
\nonumber &
B'_i:=\left\{ 
m_{r_D+(i-1)L,r_D+iL}(F'_i)<aL/2\right\}.
&\end{eqnarray}
We have used here that $\phi_z(t)\le \psi_z(t)$ (see (\ref{phi-psi})).
By the strong Markov property,
 part (a) of  Lemma \ref{psi-estimates}, and translation invariance, we see 
that for $\lambda>0$ and $n\ge 1$,

\begin{eqnarray}
\nonumber
&{\mathbb P}\left[\psi_{r_D}(F'_n,w_D)>\lambda,
 F_n< \frac{nL}{\alpha_1},
D<\infty\right]\\
\nonumber
&\!\!={\mathbb E}_w\left[{\mathbb P}_{\tau_{-r_D}w_D}\!\left[\psi_0(T_{nL})\!>\!\lambda,
 T_{nL}\!<\! \frac{nL}{\alpha_1} \right]\!
1(D\!<\!\infty)\right]\\
\label{psid2}
&\!\le\!
\frac{C_1}{\lambda}e^{-\frac{nL}{\alpha_1}\left(\alpha_1\theta
-2(\cosh\theta-1)\right)}
\end{eqnarray}
where $C_1:=3\sup_w{\mathbb E}_w[\phi_{r_D}(D),D<\infty]<\infty$ by Corollary \ref{psid}, the supremum being taken over all 
$w\in\tilde{\mathbb S}_\theta$ satisfying the conditions described, and we have
used the fact that $\psi_{r_D}(D,r_D,\eta(D))=\phi_{r_D}(D)$.  Using
(\ref{psid2})  for $n=1$, with $\lambda=p$, and Lemma \ref{est-f} we see that,

\begin{equation}
\label{extra1}
{\mathbb P}[\psi_{r_D}(T_{r_D+L})>p,D<\infty]
\le \frac{C_1}{p}e^{-\frac{L}{\alpha_1}\left(\alpha_1\theta
-2(\cosh\theta-1)\right)}
+\frac{C}{L^{M/4-1}},
\end{equation}
for some constant $C>0$.
Therefore, from (\ref{extra1}) and part $(b)$ of Lemma \ref{nol}, we have that,

\begin{eqnarray}
\nonumber
&{\mathbb P}[B_1,D<\infty]\le {\mathbb P}[\psi_{r_D}(T_{r_D+L})>p,D<\infty]
+{\mathbb P}[ m_{r_D,r_D+L}(T_{r_D+L})<aL/2]\\
\label{extra3}
&\le
\frac{C_1}{p}e^{-\frac{L}{\alpha_1}\left(\alpha_1\theta
-2(\cosh\theta-1)\right)}
+\frac{C}{L^{M/4-1}}+\frac{C}{L^{aM/(4+2M)}},
\end{eqnarray}
for some constant $C>0$.
Let us examine now the terms with $n\ge 2$ in (\ref{tail1}). Note
that in this case,
$\psi_{r_D+(n-1)L}=\psi_{r_D}+\sum_{k=1}^{n-1} \Delta_k$ where
\begin{equation} \nonumber
\Delta_k:=\psi_{r_D+kL}- \psi_{r_D+(k-1)L}
.\end{equation}
Since $\frac{1}{2^{n-1}}+\sum_{k=1}^{n-1}\frac{1}{2^{n-k}}=1$, we have,
\begin{equation}
\label{inclusion}
\left\{\psi_{r_D+(n-1)L}>p\right\}\subset
\left\{\psi_{r_D}> {p}/{2^{n-1}}\right\}
\cup\left[\cup_{k=1}^{n-1}\left\{ \Delta_k
>
{p}/{2^{n-k}}\right\}\right].
\end{equation}
Let 
\begin{equation} \nonumber
A^n_0:=\left\{\psi_{r_D}(F_n')>
{p}/{2^{n-1}}\right\},
\quad A^n_k:=\left\{
\Delta_k
( F_n')>
{p}/{2^{n-k}}\right\},
\end{equation} for $ 1\le k\le n-1$.
From
(\ref{inclusion}), for $n\ge 2$,
\begin{equation} \nonumber
B_n\subset B_{n-1}\cap \left(B'_n\cup A^n_0\cup A^n_1\cup
\cdots\cup A^n_{n-1}\right).
\end{equation}
So for $n\ge 2$,
\begin{eqnarray}\label{inc} 
 {\mathbb P}\left[B_n, D<\infty \right]
&\le& 
\sum_{k=0 }^{n-1 }{\mathbb P}\left[
A^n_k,B_{n-1},D<\infty\right]
+{\mathbb P}\left[B'_n, B_{n-1}, D<\infty\right].
\end{eqnarray}
By the strong Markov property,  we have for any $\lambda\in\mathbb R$ and  $1\le k\le n-1$,
\begin{equation}
{\mathbb P}\left[\left. F'_n-F'_k\ge  \lambda, ~D<\infty\right|{\mathcal F}_{F'_{k-1}}\right]
\le
{\mathbb P}_{w_{F'_{k-1}}}\left[ T_{(n-k+1)L}-T_L\ge  \lambda\right].
\end{equation}
Hence,  by part (c) of Lemma \ref{slowdown} and
the fact that $\alpha_1<\alpha$,
we have 
\begin{equation} \nonumber
{\mathbb P}\left[\left. F'_n-F'_k\ge   (n-k)L/\alpha_1\right|{\mathcal F}_{F'_{k-1}}\right]
\le {C}{\left((n-k)L\right)^{-M/4}}.
\end{equation}
By the strong Markov property again,
 and part (b) of  Lemma \ref{psi-estimates},
\begin{eqnarray}
\nonumber
&{\mathbb P}\left[\left.\Delta_k(F'_n)>{p}/{2^{n-k}}, F'_n-F'_k<
 {(n-k)L}/{\alpha_1}\right|{\mathcal F}_{F'_{k-1}}\right]\\
&\le 3\frac{aL }{p}2^{n-k}e^{-\frac{(n-k)L}{\alpha_1}
\left(\alpha_1\theta-2(\cosh\theta-1)\right)}.\end{eqnarray}
Therefore, for $n\ge 2$ and  $1\le k\le n-1$,
\begin{eqnarray}
\nonumber
{\mathbb P}\left[
A^n_k |{\mathcal F}_{F'_{k-1}}\right]
&\le
{\mathbb P}\left[\left.\Delta_k(F'_n)>\frac{p}{2^{n-k}}, F'_n-F'_k< \frac{1}{\alpha_1} (n-k)L\right|{\mathcal F}_{F'_{k-1}}\right]
\\
\nonumber
&+{\mathbb P}\left[\left. F'_n-F'_k\ge  \frac{1}{\alpha_1} (n-k)L\right|{\mathcal F}_{F'_{k-1}}\right]\\
\label{ank}
&\le 3\frac{aL}{p}2^{n-k}e^{-\frac{(n-k)L}{\alpha_1}
\left(\alpha_1\theta-2(\cosh\theta-1)\right)}
+\frac{C}{\left((n-k)L\right)^{M/4}}
\end{eqnarray}
From inequality (\ref{psid2}) with $\lambda=2^n/p$, Lemma \ref{est-f} and the assumption that
initially $m_{-L,0}(0)\ge aL/2$, we then obtain that for $n\ge 2$,

\begin{equation}
\nonumber
 {\mathbb P}\left[
A^n_0, D<\infty\right]
\label{a0}
\le C_1 \frac{2^n}{p}e^{-\frac{nL}{\alpha_1}\left(\alpha_1\theta
-2(\cosh\theta-1)\right)}
+\frac{C}{\left(nL\right)^{M/4-1}}.
\end{equation}
Now, for $n\ge 2$, by part $(a)$ of  Lemma \ref{nol}, the strong Markov property,
and the fact that there are $a$ particles at the rightmost visited site
at time $F'_{n-1}$,

\begin{equation}
\label{cl}
{\mathbb P}\left[\left. B'_n\right|{\mathcal F}_{F'_{n-1}}\right]
\le \frac{C}{L^{a/2} }.
\end{equation}Define a sequence 

\begin{eqnarray}
\label{asec1}
&\ &a_1:=3\frac{C}{L^{aM/(4+2M)}},\\
\label{asec2}
&\ &a_n:=\frac{4C}{((n-1)L)^{M/4-1}}\qquad {\rm for}\ n\ge 2.
\end{eqnarray} Now note that
there is a $L_0\ge C_1$, such that if
  $L\ge L_0$, for $n=1,2,\ldots$, we have that
 \begin{equation}
\label{cond1}
(3aL+C_1)\frac{2^n}{p}e^{
-\frac{nL}{\alpha_1}\left(
\alpha_1\theta-2(\cosh\theta-1)\right)}\le {C}{(nL)^{1-M/4}}\le {a_n}/{4}.
\end{equation}
(which is possible by inequality (\ref{intermediate})) and  that,  
\begin{equation} \nonumber
\sum_{n=1}^\infty a_n<1.
\end{equation} 
Let us now define $c_n:={\mathbb P}\left[B_n, D<\infty\right]$
for $n\ge 1$. We want to prove that
 the sequence $\{c_n:n\ge 1\}$ satisfies

\begin{eqnarray}
\label{cen0}
&\ &c_1\le a_1\\
\label{cen}
&\ &c_n\le a_{n}+\sum_{k=1}^{n-1}a_{n-k}c_k\qquad n\ge 2.
\end{eqnarray}
From (\ref{extra1}), (\ref{cond1}) and the fact that
$1/L^{M/4-1}\le 1/L^{aM/(4+2M)}$ (which follows from
(\ref{mcond})) 
 note that (\ref{cen0}) is satisfied.
Now note that by inequality (\ref{cl}),  whenever $L\ge L_0$, for
$n\ge 2$ we have that
\begin{equation}
\label{aa1}
{\mathbb P}[B'_n, B_{n-1}, D<\infty]\le \frac{C}{L^{a/2}}
{\mathbb P}[B_{n-1}, D<\infty]\le
 a_1{\mathbb P}[B_{n-1}, D<\infty].
\end{equation}
 Inequality (\ref{a0}) and condition
(\ref{cond1}) imply ${\mathbb P}\left[
A^n_0,D<\infty\right]
\le {a_{n}}/{2}$ and
${\mathbb P}\left[
A^{n}_1,D<\infty\right]
\le {a_{n}}/{2}$. Hence, for $n\ge 2$,

\begin{equation}
\label{aa01}
{\mathbb P}[A^n_0,B_{n-1},D<\infty]
+{\mathbb P}[A^n_1,B_{n-1},D<\infty]\le a_n.
\end{equation}
Similarly for $2\le k\le n-2$ we have  ${\mathbb P}\left[
A^n_{k}|{\mathcal F}_{F'_{k-1}}\right]
\le a_{n-k+1}$. Thus, since $B_{n-1}\subset B_{k-1}$
 for $2\le k\le n-2$,

\begin{equation}
\label{aak}
{\mathbb P}[A_{k}^n, B_{n-1}, D<\infty]\le
{\mathbb P}[A_{k}^n, B_{k-1}, D<\infty]\le a_{n-k+1}
{\mathbb P}[B_{k-1},D<\infty].
\end{equation}
Also, by inequality (\ref{ank}) and condition
(\ref{cond1}) for $n\ge 2$, we have  ${\mathbb P}\left[
A^n_{n-1}|{\mathcal F}_{F'_{n-2}}\right]
\le \frac{a_{2}}{2}$. Thus, since $B_{n-1}\subset B_{n-2}$, for $n\ge 3$,

\begin{equation}
\label{aa2}
{\mathbb P}[A_{n-1}^n, B_{n-1}, D<\infty]\le {\mathbb P}[A_{n-1}^n, B_{n-2}, D<\infty]\le a_2
{\mathbb P}[B_{n-2},D<\infty].
\end{equation}
For $n=2$, (\ref{cen}) now follows after substituting estimates (\ref{aa1})
and (\ref{aa01})  in inequality (\ref{inc}), for $n=3$
after substituting (\ref{aa1}), (\ref{aa01}) and (\ref{aa2}) in
inequality (\ref{inc}) while
for $n\ge 4$ it follows after substituting (\ref{aa1}), (\ref{aa01}),
(\ref{aak}) and (\ref{aa2})
 in inequality (\ref{inc}).

But remark that,

\begin{equation} 
\label{asum}
\sum_{n=1}^\infty n^{M/4-3}a_n<\infty.
\end{equation}
Hence, by Lemma \ref{series}, (\ref{cen0}), (\ref{cen}) and inequality (\ref{tail1}) we conclude that,

\begin{equation} \nonumber
\sum_{n=1}^\infty n^{M/4-3}{\mathbb P}\left[J_{r_D}>n,
D<\infty\right]<\infty.
\end{equation}
This implies that $\limsup_{n\to\infty}n^{M/4-3}{\mathbb P}\left[J_{r_D}>n,
D<\infty\right]=0$. Thus, there exists a constant $C>0$, such that for
every $n\ge 1$ it is true that
${\mathbb P}\left[J_{r_D}>n,
D<\infty\right] \le C/n^{M/4-3}$.
This together with the monotonicity in $t$ of the
expression ${\mathbb P}\left[J_{r_D}>t,
D<\infty\right]$, finishes the proof.
\smallskip

\noindent {\it Proof of part $(b)$}. This time, in analogy
with (\ref{tail1}), note that for  $n= 1,2,\ldots$,
\begin{equation}
\label{tail2}
{\mathbb P}[J_0>n, U=\infty]
\le {\mathbb P}\left[ B_n, U=\infty\right],
\end{equation}
where again we have dropped the subscript on $\mathbb P_w$ but now
\begin{eqnarray}&
\nonumber
B_n:=\cap_{i=1}^n \left\{\psi_{(i-1)L}(T_{iL},w)>p\right\}\cup B'_i,&\\
\nonumber &
B'_i:=\left\{ 
m_{(i-1)L,iL}(T_{iL})<aL/2\right\}.
&\end{eqnarray}
An analysis similar to that of part $(a)$ proves (\ref{cen0}),
(\ref{cen}) and (\ref{asum}) with $c_n:={\mathbb P}\left[ B_n, U=\infty\right]$
for $n\ge 1$ and $\{a_n:n\ge 1\}$ as in (\ref{asec1}) and (\ref{asec2}).
Part $(b)$ now follows by Lemma \ref{series} as in part $(a)$.

\smallskip

\noindent {\it Proof of (\ref{partc})}. Note that
$\{J_0=\infty\}\subset\{J_{L}=\infty\}$. Now,
for every $n\ge 1$,
$$
{\mathbb P}_w[J_L=\infty]\le
{\mathbb P}_w[J_L=\infty, m_{-n,L}(T_L)\ge aL/2]+
{\mathbb P}_w[m_{-n,L}(T_L)< aL/2].
$$
Now, following the proof of part $(a)$, it is possible to
show that ${\mathbb P}_w[J_L>t,m_{-n,L}(T_L)\ge aL/2]\le C 
t^{3-M/4}$, for some constant $C>0$.  Hence, for every $n\ge 1$,

$$
{\mathbb P}_w[J_L=\infty]\le
{\mathbb P}_w[m_{-n,L}(T_L)< aL/2].
$$
Taking the limit as $n\to\infty$ and using
part $(c)$ of Lemma \ref{nol}
we finish the proof.
\end{proof}

\smallskip

\smallskip

\begin{corollary}
\label{sfinite} For every nontrivial initial condition $w\in\tilde {\mathbb S}_\theta$, it is true that,

\begin{equation}
\label{sf1}
S_1<\infty\qquad  {\mathbb P}_w-{\rm a.s.},
\end{equation}
and for every $k\ge 2$,

\begin{equation}
\label{sf2}
{\mathbb P}_w[D_{k-1}<\infty, S_k<\infty]={\mathbb P}_w[D_{k-1}<\infty].
\end{equation}
\end{corollary}

\begin{proof} Assertion (\ref{sf1}) is a consequence of
(\ref{partc}) of Lemma \ref{fundamental} and the
fact that $\{S_1<\infty\}=\{J_0<\infty\}$. Similarly, 
assertion (\ref{sf2}) follows directly from part $(a)$
of Lemma \ref{fundamental} and the fact that
$\{D_{k-1}<\infty, S_k<\infty\}=\{D_{k-1}<\infty, J_{r_{D_{k-1}}}
<\infty\}$.

\end{proof}
\smallskip

\subsection{Variance bounds for the regeneration times and positions} 
\label{regen-tp}
In this subsection we will prove 
 Proposition \ref{moment2}. Let us first prove assertion (\ref{weight1})
of Proposition \ref{moment2}.
By Corollary  \ref{sfinite}, note that for every $k\ge 1$,

$$
{\mathbb P}_w[\kappa=\infty]\le {\mathbb P}_w[D_k<\infty].
$$
But by the strong Markov property and Lemma \ref{uniform},
the righthand side of the above inequality is bounded by $(1-\delta)^k$.
It follows that,

$$
{\mathbb P}_w[\kappa=\infty]\le (1-\delta)^k,
$$
for every $k\ge 1$. Taking the limit when $k$ tends to infinity concludes
the proof of (\ref{weight1}) of Proposition \ref{moment2}.

To prove (\ref{moments-delta}) we will need the following lemma.

\smallskip
\begin{lemma}

\label{moment} For every $\epsilon>0$, there is a
constant $C$,  $0<C<\infty$, such that

\begin{equation}
\label{kappamoment}
{\mathbb P}_{a\delta_0}
\left[\left. \kappa>t\right|U=\infty\right]\le Ct^{-M/4+3+\epsilon}.
\end{equation}

\end{lemma}

\begin{proof}[proof] Without loss of generality we will assume that
initially,

\begin{equation}
\label{r-in}
r_0=0.
\end{equation}
By the fact that $\kappa<\infty$, a.s., we can write,

\begin{equation} \nonumber
{\mathbb P}\left[\left. \kappa>t \right| U=\infty\right]
=\sum_{k=1}^\infty {\mathbb P}\left[\left. S_k >t, K=k \right| U=\infty\right],
\end{equation}
where we have dropped the subscript $a\delta_0$ in
 ${\mathbb P}_{a\delta_0}$.
Applying recursively the strong Markov property
to the stopping times $\{S_j:j\ge 1\}$ we see that for every $k\ge 1$,

\begin{equation} \nonumber
{\mathbb P}\left[\left. S_k >t, K=k \right| U=\infty\right]
\le (1-\delta)^{k-1},
\end{equation}
where $\delta>0$ is given by
 Lemma \ref{uniform}. Let $0<\beta<1/2$. For any $l>0$ we therefore have,

\begin{equation}
\label{decomp}
{\mathbb P}\left[\left. \kappa>t \right| U=\infty\right]
\le
\sum_{k=1}^l {\mathbb P}\left[\left. t<S_k <\infty \right| U=\infty\right]
+\delta^{-1}(1-\delta)^l.
\end{equation}
Let $0<\gamma<1$ and consider the event,

\begin{equation} \nonumber
A_k:=\{r_{D_1}-r_{S_1}<t^\gamma, r_{D_2}-r_{S_2}<t^\gamma,\ldots,
r_{D_{k-1}}-r_{S_{k-1}}<t^\gamma, S_k<\infty\}
\end{equation}
On $A_k$ we have,
$
r_{S_k}\le k t^\gamma+L\sum_{j=0}^{k-1} J_{r_{D_j}}$, where
we adopt the convention $D_0:=0$, so that $r_{D_0}=0$ by (\ref{r-in}).
Since $\tilde r_t\le r_t$, if $U=\infty$, then $r_t\ge \lfloor \alpha_2 t\rfloor$
for all $t>0$. Therefore, on $A_k\cap\{U=\infty\}$,

\begin{equation} \nonumber
\lfloor \alpha_2 S_k\rfloor\le k t^\gamma+L\sum_{j=0}^{k-1} J_{r_{D_j}}.
\end{equation}
Now define, the event

\begin{equation}
B_k:=\{J_{r_{D_0}}
<t^\gamma, J_{r_{D_1}}<t^\gamma,\ldots,
J_{r_{D_{k-1}}}<t^\gamma,S_k<\infty\}.
\end{equation}
Then on $A_k\cap B_k\cap\{U=\infty\}$ we have,

\begin{equation} \nonumber
\lfloor \alpha_2 S_k\rfloor\le k t^\gamma(1+L).
\end{equation}
Hence for $t>(lt^\gamma (1+L)+1)/\alpha_2$ and $k\le l$,

\begin{equation} \nonumber
{\mathbb P}[t<S_k<\infty, A_k, B_k|U=\infty]=0
\end{equation}
and therefore,

\begin{equation}
\label{lemmas}
{\mathbb P}[t<S_k<\infty|U=\infty]
\le {\mathbb P}[A_k^c,S_k<\infty|U=\infty]
+
{\mathbb P}[B_k^c,S_k<\infty|U=\infty].
\end{equation}
Using part $(a)$ of  Lemma \ref{fundamental} to bound the
probability of the event $\{J_{r_0}\ge t^\gamma\}=
\{J_0\ge t^\gamma\}$ and part $(b)$ to
bound the probability of the events $\{J_{r_{D_j}}\ge t^\gamma\}$,
for $1\le j\le k-1$, we can see that the second term of the righthand side of
inequality (\ref{lemmas}) is bounded by
$C k t^{-\gamma (M/4-3)}$. On the other hand, by Lemma \ref{rd}, the first term
is bounded by $C k t^{-\gamma M/2}$. Choosing $l=C_1\log t$
with $C_1=\left(M/4-3\right)(\log (1-\delta)^{-1})^{-1}$ and
$\gamma$ close enough to $1$ we obtain (\ref{kappamoment}).
\end{proof}
\smallskip

 {\em Proof 
 of  (\ref{moments-delta}) of Proposition \ref{moment2}.} 
The assertion for $\kappa$ of
 (\ref{moments-delta}) follows from Lemma \ref{moment}
 noting that $M\ge 21$ (by condition
(\ref{mcond})) and that for $r_{\kappa}$  from 
Lemma \ref{moment} and  (\ref{estimater}).

\section{Construction and Feller property}
\label{fellerp}
Throughout, $\theta>0$ is arbitrary, $P$ is the joint law of the
independent random walk used to define the process for finite
initial conditions and $E$ the corresponding expectation.
By our construction note that  $r^\ell_t$ is increasing
in $\ell$ and hence we can define
\begin{equation} \nonumber
r_t:=\lim_{\ell\to\infty} r^\ell_t.
\end{equation}
We will see that for every $t\ge 0$, a.s. $r_t<\infty$.

 Consider $f_\theta(\eta^\ell)$ where
\begin{equation} \nonumber
f_\theta(\eta)= \sum_{x} \eta(x) e^{ \theta x }.
\end{equation}
We compute
\begin{equation} \nonumber
{\mathcal L}f_\theta(\eta^\ell)=   \sum_{x} \eta^\ell(x) 
 e^{\theta x } \left[ e^\theta + e^{-\theta} -2 + ((a+1)e^\theta-1) 1(x=r)\right].
 \end{equation}
Hence if we let $\lambda_{1,\theta} = e^\theta + e^{-\theta} -2$ and
$\lambda_{2,\theta} = (a+1)e^\theta + e^{-\theta} -2$ then
\begin{equation} \nonumber
\lambda_{1,\theta} f_\theta\le {\mathcal L}f_\theta \le \lambda_{2,\theta} f_\theta.
\end{equation}
In particular, 
\begin{equation}\label{11}
E[ f_\theta(\eta^\ell(t))~|~\mathcal F_0] \le e^{\lambda_{2,\theta}t} 
f_\theta(\eta^\ell(0)).
\end{equation}
In addition, $ f_\theta(\eta^\ell(t))$ is 
a nonnegative sub-martingale and therefore by Doob's inequality,
\begin{equation}\label{12}
P(\sup_{0\le s\le t}  f_\theta(\eta^\ell(s)) \ge e^{\gamma \theta t}~|~\mathcal F_0)
\le e^{-\gamma\theta t} E[  f_\theta(\eta^\ell(t))~|~\mathcal F_0].
\end{equation}
Since $r^\ell_t$ is the rightmost site which has been occupied up to
time $t$ we have $\sup_{0\le s\le t} f_\theta(\eta^\ell(s)) \ge e^{\theta r^\ell_t}$.  Hence from (\ref{11}) and (\ref{12}) we have
\begin{equation}\label{estimate-r}
P( r^\ell_t \ge \gamma t~|~\mathcal F_0) \le e^{ -c_{\gamma,\theta}
 t}f_\theta(\eta^\ell(0))
\end{equation}
where $c_{\gamma,\theta}=\gamma \theta-\lambda_{2,\theta}$.
This proves that for each $\ell$ and $t\ge 0$, a.s. $r^\ell_t<\infty$
and hence $\lim_{n\to\infty}\tau_n=\infty$.
Also, taking the limit when $\ell\to\infty$ in (\ref{estimate-r}), we obtain,
\begin{equation}\label{estimate-r2}
P( r_t \ge \gamma t~|~\mathcal F_0) \le e^{ -c_{\gamma,\theta}
 t}f_\theta(\eta(0))
\end{equation} 
This proves Lemma \ref{front1} of Section \ref{boundsd}. 
Furthermore, if $(r,\eta)\in \mathbb S'_\theta$ then $f_\theta(\eta)<\infty$
so we have $
r_t<\infty$  a.s.

  Choose now $\gamma$ large
enough so that we have $c_{\gamma,\theta} >0$.
Define for each $y= 1,2,\ldots$, 
\begin{equation}
T_y:=\inf\{t\ge 0:r_t=y\}.
\end{equation} 
We have $r_{T_y\land t}=
\lim_{\ell\to\infty} r^\ell_{T_y\land t}$.  Let $\ell_k$ be the
smallest natural number such that $r_{T_k\land t}=r^\ell_{T_k\land t}$
for all $\ell\ge \ell_k$. Then if $\bar\ell=\max\{\ell_1,\ldots, \ell_{r_t}\}$,
the front $r_{\cdot\land t}$ generated by the initial condition $\eta^\ell$
up to time $t$, does not depend on $\ell$ if $\ell\ge \bar\ell$. This means that
particles that are initially at any site $x\le r-\bar\ell$, never visit
any site to the right of the front before time $t$. Using 
attractiveness, we can then conclude that the
sequence $\eta^\ell(s)$ is increasing for $s\le t$ and $\ell\ge \bar\ell$. 
Therefore, 
\begin{equation} \nonumber
\eta(t):=\lim_{\ell\to\infty}\eta^\ell(t),
\end{equation}
exists almost surely. 
Taking the limit when $\ell\to\infty$ in (\ref{11}) and using
Fatou's Lemma we see that,
\begin{equation}\label{1122}
E[ f_\theta(\eta(t))~|~\mathcal F_0] \le e^{\lambda_{2,\theta}t} 
f_\theta(\eta(0)).
\end{equation}
Noting that $r_t$ is increasing, this shows that $(r_t,\eta(t))$ stays in $\mathbb S'_\theta$.  Hence we have shown that
$(r_t,\eta(t))$ is a Markov process on $\mathbb S'_\theta$.

We next want to show that it satisfies the Feller property.  We
need some more preliminary estimates.
Note that
 \begin{eqnarray}
 \nonumber
&&{\mathcal L}f_\theta^2-2f_\theta {\mathcal L}f_\theta 
\\ 
\nonumber
 && =  \sum_{x} \eta(x) 
 e^{2\theta x } \left[ (e^\theta-1)^2 + (e^{-\theta}-1)^2 +
 ( (a^2-1)e^{2\theta}  -2(a-1) e^\theta ) \right]\\
&& \le  \lambda_{3,\theta} f_{2\theta}
 \end{eqnarray} for some $\lambda_{3,\theta} <\infty$.
 Hence if $M_\theta(t) = f_\theta(\eta^\ell(t)) - \int_0^t
 {\mathcal L}f_\theta( \eta^\ell(s))ds $ then
$\left\langle M_\theta(t)\right\rangle = \int_0^t\left(
{\mathcal L} f^2_\theta-
2f_\theta {\mathcal L}f_\theta\right)(\eta^\ell(s))ds $, and then,
 \begin{equation} \nonumber
 E[ M_\theta^2(t) ~|~\mathcal F_0] 
 \le 2\lambda_{3,\theta} \int_0^t E[ f_{2\theta}(\eta^\ell(s))~|~\mathcal F_0]
 ds \le 2\lambda_{3,\theta}\lambda_{2,2\theta}^{-1} e^{\lambda_{2,2\theta} t } f_{\theta}(\eta^\ell(0)). 
 \end{equation}
 Here we used in the last inequality that $f_{2\theta}(\eta^\ell(0))\le 
 f_{\theta}(\eta^\ell(0))$, since there are no particles initially to the
 right of the origin.  In particular, by Chebyshev's inequality, 
 \begin{equation}\label{estimate-f}
 P(f_\theta(\eta^\ell(t)) < e^{\lambda_{1,\theta} t} f_\theta(\eta^\ell(0)) 
 -A~|~\mathcal F_0) \le A^{-2} e^{\lambda_{4,\theta} t} f_\theta(\eta^\ell(0))
 \end{equation}
  for some $\lambda_{4,\theta}<\infty$. We can pass to  the limit to obtain,
\begin{equation}\label{estimate-f2}
 P(f_\theta(\eta(t)) < e^{\lambda_{1,\theta} t} f_\theta(\eta(0)) 
 -A~|~\mathcal F_0) \le A^{-2} e^{\lambda_{4,\theta} t} f_\theta(\eta(0)).
 \end{equation}
This proves that we have a well-defined Markov  process  starting from any initial data
in $\mathbb S'_\theta$.  Next we show the process satisfies the Feller property.  We start
by identifying the compact sets of $\mathbb S'_\theta$.

\begin{lemma}
The compact sets $K$ in $\mathbb S'_\theta$ are those which are 
closed, bounded in norm $\|(r,\eta)\| = d((r,\eta),(0,0))$,
and have {\em uniform tails} i.e. 
\begin{equation}
\lim_{N\to\infty} \sup_{(r,\eta)\in K} \sum_{x\le r-N}
e^{\theta (x-r)} \eta(x) =0.
\end{equation}\label{a}
\end{lemma}
\begin{proof}
Suppose $(r,\eta)_i$, $i=1,2,\ldots$ are elements of such a $K$.
Since $K$ is bounded, 
we can find a weakly convergent subsequence.    We have to show they converge in norm as well.  Relabeling, we can 
call the weakly convergent subsequence $(r,\eta)_i\to (r,\eta)$.  Also,
note that $r_i=r$ for large enough $i$.  So without loss
of generality assume $r_i=r$.
Let $\epsilon>0$ and choose $N_0$ such that 
\begin{eqnarray}&
\sup_{i} \sum_{x<r-N_0} e^{\theta(x-r)} \eta_i(x) <\epsilon. 
&
\end{eqnarray}
Now choose $i_0$ so that for $i\ge i_0$,
\begin{eqnarray}
&
\sum_{r-N_0\le x\le r} e^{\theta(x-r)} |\eta_i(x) -\eta(x)|<\epsilon.
&
\end{eqnarray}
It follows that for $i\ge i_0$,
\begin{equation} \nonumber
d((r,\eta)_i, (r,\eta)) < 3\epsilon.
\end{equation}
Since $K$ is closed $(r,\eta)\in K$ and  such a set is compact.  
Suppose on the other hand that a subset $K$ of $\mathbb S'_\theta$ is
compact.  Fix $\epsilon>0$.  Let
\begin{equation} \nonumber
B_N = \{ (r,\eta)  : \sum_{x\le r-N}e^{\theta(x-r)} \eta(x) <\epsilon\}.
\end{equation}
$B_N$ are open sets whose union is $\mathbb S'_\theta$.  So $B_N$
are an open cover of $K$.  Since $K$ is compact, there is 
a finite sub-cover, and hence an $N$ such that
$K\subset B_N$.  In other words, $K$ has uniform tails.
\end{proof}

\begin{lemma}
For each $t>0$, $\epsilon>0$ and $K\subset \mathbb S'_\theta$ compact there exists
a compact $K_0\subset \mathbb S'_\theta$ such that 
\begin{equation}
P( (\eta(t), r_t) \in K~|~ (\eta(0), r_0=0)\not\in K_0)<\epsilon.
\end{equation}\label{b}
\end{lemma}
\begin{proof} Fix $t>0$, $\epsilon>0$ and $K\subset \mathbb S'_\theta$ compact. 
By Lemma \ref{a}, $K$ is of the
form 
\begin{eqnarray}
&K=\{ (r,\eta)~|~ \sum_{x\le r-N(m)}  e^{\theta(x-r)}\eta(x) \le  B/m, {\rm ~for ~all~} m=1,2,\ldots\}&
\end{eqnarray}
 for some
$B<\infty$ and some $N(m)\uparrow \infty$ as $m\to\infty$ with $N(1)=0$.  
We have to find a $B_0$ and $N_0(\cdot)$ so that
for each $m_0=1,2,\ldots$, if $\sum_{x\le -N_0(m_0)} e^{\theta x} \eta(0,x) >B_0/m_0$ then $P(\sum_{x\le r_t-N(m)}  e^{\theta(x-r_t)}\eta(t,x) \le  B/m {\rm ~for ~all~} m=1,2,\ldots\})<\epsilon$.

Choose $\gamma$ large enough such that 
\begin{eqnarray}\label{gamma}
& P(r_t > \gamma t) <\epsilon/2.&
\end{eqnarray}

We start with $m_0=1$.  Use  (\ref{estimate-f2}) with  $A^2= f_\theta(\eta(0)) e^{\lambda_{4,\theta} t}\epsilon/2$.  We can find $B_0$ so
that  if $f_\theta(\eta(0)) > B_0$
then $e^{\lambda_{1,\theta} t} f_\theta(\eta(0)) -A> B e^{\gamma t}$ and
hence from (\ref{estimate-f2}) and (\ref{gamma}), if $\sum_{x\le 0} e^{\theta x} \eta(0,x) >B_0$,
\begin{eqnarray} \label{112}
& P(\sum_{x\le r_t} e^{\theta (x-r_t) } \eta(t,x)  \le B) <\epsilon. &
\end{eqnarray}  

Next we consider the case $m_0>1$.  It is not hard to see that for each $N$, there exists $A=A(t)<\infty$ such that
\begin{eqnarray}\label{mzero}
& P(\sum_{x\le -N} e^{\theta x} \eta(t,x) \le \frac12 \sum_{x\le -N-A} 
e^{\theta x} \eta(0,x)) < \epsilon/2. &
\end{eqnarray}
Indeed, the left hand side of the event in (\ref{mzero})  is only smaller if we
suppress the branching.  If we temporarily denote $x_i(0)$ the initial positions of particles to the left of $-N -A$ then we have continuous time random
walks and the event is that
$\sum_i e^{\theta x_i(t)} 1(x_i(t) \le -N) \le \frac12
\sum_i e^{\theta x_i(0)} $.  We can assume that $\sum_{x\le 0} e^{\theta x}\eta(0,x) \le B_0$,
for otherwise we have (\ref{112}).  Then it is clear that there exists
an $A$ such that $P(\cup_i \{x_i(t) >-N\}) <\epsilon/4$.
Hence we only need to show that 
$P(\sum_i e^{\theta x_i(t)} \le \frac12
\sum_i e^{\theta x_i(0)})<\epsilon/4$  which is easy to deduce from the fact that $e^{ - \lambda_{1,\theta} t }\sum_i e^{ \theta x_i(t)}$ are 
martingales.  This proves (\ref{mzero}).

From (\ref{gamma}) and (\ref{mzero}) with 
 $N_0(m) = N( \lfloor B B_0^{-1} 2me^{\gamma t}\rfloor +1 )+A$,
 we have that if $\sum_{x\le - N_0(m)} e^{\theta x} \eta(0,x) > B_0/m$ and
 $\sum_{x\le 0} e^{\theta x}\eta(0,x) \le B_0$ then
\begin{eqnarray} &
P(\sum_{x\le r_t- N(m')} e^{\theta (x-r_t)} \eta(t,x) \le B/m')
\le \epsilon,&\end{eqnarray}
for $m'= \lfloor B B_0^{-1} 2me^{\gamma t}\rfloor +1$.
This completes the proof. 
\end{proof}

\begin{lemma}
\label{cont}
 For each $\epsilon>0$ and $t>0$ there exists $\delta>0$ such 
that if $(r,\eta)$ and $(r,\eta')$ are any 
two configurations of particles on $\mathbb S'_\theta$ with
$\sum_{x\le r} e^{\theta x}|\eta(x)-\eta'(x)| <\delta$, there is stopping time $\tau$ and a coupling of two copies $(r_s, \eta(s))$ and $(r'_s,\eta'(s))$ of our Markov process with generator
${\mathcal L}$ for 
$0\le s\le \tau$ satisfying  
\begin{enumerate}
\item  $P(\tau < t) <\epsilon$.
\item $E[ d( (r_t,\eta(t)), (r'_t,\eta'(t))) \mathbf 1_{\tau>t}] <\epsilon$.
\item $P[r_0=r'_0=r,\eta(0)=\eta,\eta'(0)=\eta']=1$. 
\end{enumerate}
Here $P$ is the coupling
measure and $E$ the corresponding expectation.
\end{lemma}
\begin{proof}
Consider the  difference
$\zeta= \eta
-\eta'$.   We have $\sum_{x\le 0} e^{\theta x} |\zeta(x)| <\delta$ so choosing $\delta$ sufficiently small we have $\zeta(x)=0$ for $x\in \{-L,\ldots,0\}$ for some large $L$.  We attempt to couple the two processes by moving the particles together whenever possible.  Then positive and negative parts of $\zeta$
move as independent random walks of positive and negative type, the
two types annihilating on contact and the coupling succeeds up to the first time
$\tau$ when a particle of either type hits $r_s$.   It is easy to choose
$\delta$ small enough, and therefore $L$ large enough, so that 
(1) is satisfied.  To prove (2), note that up to time $\tau$, 
$d((r_t, \eta(t)), (r'_t,\eta'(t))) = \sum_{x\le r_t} e^{\theta(x-r_t)} |\zeta(t,x)|$ and
we can get an easy upper bound by using $r_t\ge 0$ and letting $\bar\zeta(t)$ be the process obtained by starting with $|\zeta(0)|$ and using the 
same random walks, but dropping the signs and the annihilations. 
 Then  $\sum_{x} e^{\theta (x-r_t)} |\zeta(t,x)|\le 
\sum_x e^{\theta x} \bar\zeta(t,x)$.
 (2) follows since $e^{-\lambda_{1,\theta}t}\sum_{x}
 e^{\theta x} \bar\zeta(t,x)$ is a martingale. \end{proof}

\begin{proposition} Let $(r,\eta)\in{\mathbb S}'_\theta$ and $P_{r,\eta}$
the law of the process $\{(r_t,\eta(t):t\ge 0\}$ with initial condition
$(r,\eta)$ under $P$. Then,
$P_t g( r,\eta) = E_{r,\eta}[ g(r_t,\eta(t))] $, $t\ge 0$ form a Feller semi-group on $\mathbb S'_\theta$, where $E_{r,\eta}$ is the expectation associated
to $P_{r,\eta}$.
\end{proposition}
\begin{proof}
 Suppose that $g$ is continuous and vanishes at infinity and let $\epsilon>0$.  
 In particular $|g|\le B<\infty$.
 There is a compact $K$ such that 
 $|g(r,\eta)|<\epsilon/2$ for $(r,\eta) $ in the complement of
 $K$.  By Lemma \ref{b} there is a compact $K_0$ such that
 $P((r_t,\eta(t))\in K~|~(r_0,\eta(0))\not\in K_0)
 <\epsilon/2B$.  So if $(r,\eta)\not\in K_0$,
 \begin{equation} \nonumber
 P_tg( r,\eta) = E_{r,\eta}[ g(r_t,\eta(t)) \mathbf 1_{(r_t,\eta(t))\in K}]+E_{r,\eta}[ g(r_t,\eta(t)) \mathbf 1_{(r_t,\eta(t))\not\in K}] <\epsilon.
 \end{equation}
  This proves that $P_tg$
vanishes at infinity as well.  

Next we show that $P_tg$ is continuous.  Since $g$ is
continuous and vanishes at infinity, it is uniformly continuous.  So we
can choose $\epsilon_0>0$ so that $d((r,\eta),(r',\eta'))<\epsilon_0$ implies
$|g(r,\eta)-g(r',\eta')|<\epsilon/3.$   By Lemma \ref{cont}, there exists
$\delta$ such that if $d((r,\eta),(r',\eta'))<\delta$, then
there is a stopping time $\tau\ge 0$ such that $P(d((r_t,\eta(t)),(r'_t,\eta'(t)))>\epsilon_0,\tau>t)<\epsilon/(6B)$ and $P(\tau<t)\le \epsilon/(3B)$.  Hence 
$|E_{r,\eta} [ g(r_t,\eta(t))] - E_{r',\eta'}[g(r'_t,\eta'(t))]|
\le E[ |g(r_t,\eta(t))-g(r'_t,\eta'(t))|{\mathbf 1}_{ \tau>t}] +
2BP(\tau<t) \le \epsilon$.  This proves that $P_tg$ is continuous
as well.  
\end{proof}

\end{document}